\documentclass[11pt,A4,twoside]{article}\usepackage[margin=1.4in]{geometry}

\usepackage{cancel,multicol}
\usepackage{enumerate}
\usepackage{makeidx}

\usepackage{savesym,amsmath,ulem}\def\em{\it}
\savesymbol{subequations}
\usepackage{cases}
\usepackage{ifthen,bm}
\usepackage{dsfont}

\usepackage{pgf}
\usepackage{latexsym, amsfonts, amssymb}
\usepackage{alltt,latexsym, verbatim}
\usepackage{float}
\usepackage{amsfonts}
\usepackage{xspace}

\usepackage{natbib}
\def\citeN{\citet*}\def\shortciteN{\citet*}

\newtheorem{theorem}{Theorem}[section]
\newtheorem{lem}{Lemma}[section]
\newtheorem{pro}{Proposition}[section]
\newtheorem{cor}{Corollary}[section]
\newtheorem{conj}{Conjecture}[section]

\newtheorem{rem}{Remark}[section]
\newtheorem{com}{Comments}[section]
\newtheorem{ex}{Example}[section]
\newtheorem{defi}{Definition}[section]
\newtheorem{hyp}{Assumption}[section]

\numberwithin{equation}{section}

\newcommand{\bt}{\begin{theorem}}\newcommand{\et}{\end{theorem}}
\newcommand{\bl}{\begin{lem}}\newcommand{\el}{\end{lem}}
\newcommand{\bp}{\begin{pro}}\newcommand{\ep}{\end{pro}}
\newcommand{\bcor}{\begin{cor}}\newcommand{\ecor}{\end{cor}}
\newcommand{\bconj}{\begin{conj}}\newcommand{\econj}{\end{conj}}

\newcommand{\bd}{\begin{defi} \rm }\newcommand{\ed}{\end{defi} }
\newcommand{\brem }{\begin{rem} \rm }\newcommand{\erem }{\end{rem}}
\newcommand{\bcom}{\begin{com} \rm }\newcommand{\ecom }{\end{com}}
\newcommand{\brems }{\begin{rem} \rm }\newcommand{\erems }{\end{rem}}
\newcommand{\bex}{\begin{ex} \rm }\newcommand{\eex}{\end{ex}}
\newcommand{\bhyp}{\begin{hyp} \rm }\newcommand{\ehyp}{\end{hyp}}

\def\proof{\noindent \textbf{\emph{\textbf{Proof}.$\qqq$}}}
\def\finproof {\hfill $\Box$ \vskip 5 pt }

\def \Int{\displaystyle\int}

\def \be{\begin{eqnarray}}
\def \ee{\end{eqnarray}}
\def \b*{\begin{eqnarray*}}
\def \e*{\end{eqnarray*}}

\def \[{[\,\!\![}
\def \]{]\,\!\!]}

\def \1{{\bf 1}}

\def \proof{{\noindent \bf Proof. }}

\def\F{{\cal T}}\def\F{{\mathfrak{G}}}\def\F{{\mathfrak{F}}}

\newcommand{\beqa}{\begin{eqnarray}}
\newcommand{\eeqa}{\end{eqnarray}}

\def\N{N}\def\N{{\mathbb N}}
\def\R{{\mathbb R}}
\def\proof{\noindent {\it Proof. $\, $}}
\def\finproof {\hfill $\Box$ \vskip 5 pt }

\def\I{\mathds{1}}
\def\sp{\,,\ \,}

\def\bal{\begin{aligned}}
\def\eal{\end{aligned}}

\def\cL{\mathcal{L}}

\def\hat{\widehat}

\def\mym{m}

\def\xitau2{\xi_{(\tauu)}}\def\xitau2{\xi}
\def\xiitau2{\xi^i_{(\tauu)}}\def\xiitau2{\xi^i}
\def\xintau2{\tilde{\xi}_{(\tauu)}}\def\xintau2{\tilde{\xi}}

\def\Ltau2{\tilde{\xi}(\theta_{\mym},\theta_1)}
\def\chiitau2{P^i_{\theta}}

\newcommand{\beq}{\begin{eqnarray*}}
\newcommand{\eeq}{\end{eqnarray*}}

\def\mynu{\nu}
\def\mylambda{\lambda}
\def\mylambda{\lambda}

\def\mym{-}

\def\theh{h}

\def\tilde{\widetilde}
\def\cadlag{c\`adl\`ag }

\def\bX{b}

\def\sigmaX{\sigma}
\def\emph{}

\def\gg{{\mathbb G}}\def\gg{{\mathbb F}}
\def\G{{\mathfrak{G}}}\def\G{{\mathfrak{F}}}
\def\ff{\ff}

\def\ff{\ff^*}
\def\F{{\mathfrak{F}}^*}

\def\gg{{\mathfrak{G}}}\def\gg{\gg}
\def\G{{\mathfrak{G}}}
\def\ff{{\mathfrak{F}}}\def\ff{\ff}
\def\F{{\mathfrak{F}}}

\def\E{\mathbb{E}}

\def\Rf{\mathfrak{R}}\def\Rf{\mathfrak{r}}
\def\Rf{\mathfrak{R}}\def\Rf{\mathfrak{r}}\def\Rf{\mathfrak{r}^b}

\def\thisR{R}\def\thisR{\rho}\def\thisR{R}\def\thisR{\rho^b}

\def\tauu{\vartheta}
\def\bthet a{\tau}
\def\tb{{\bar{\theta}}}\def\tb{\tauu}

\def\tb{\tauu}\def\tb{{\bar{\theta}}}

\def\ind{\mathds{1}}
\def\e{z}\def\e{e}

\def\theU{U}

\def\theN{\mathbb{N}_n}\def\theN{N}

\def\Ep{{E^{\tp}}}\def\Ep{{E^{b}}}\def\Ep{\Et}

\def\mylambda {\gamma}

\def\ind {1\!\!1}\def\ind{\mathds{1}}
\def\I{\ind}

\def\F{{\mathfrak{F}}}

\def\G{{\mathfrak{G}}}

\def\ff{{\mathbb F}}

\def\gg{{\mathbb G}}

\def\P{\mathbb P}
 \def\Q{\mathbb Q}

\def\E {{\mathbb E} }

\def\N{{\mathbb N}}

\def\finproof {\hfill $\Box$ \vskip 5 pt }

\def\bal{\begin{aligned}}
\def\eal{\end{aligned}}

\def\finproof {\hfill $\Box$ \vskip 5 pt }

\def\Proba{\Proba}\def\Proba{\mathbb{Q}}

\def\sp{,\ \, }
\def\r#1{(\ref{#1})}

\def\xitau2{\xi_{(\tauu)}}
\def\xiitau2{\xi^i_{(\tauu)}}
\def\xintau2{\tilde{\xi}_{(\tauu)}}
\def\Ltau2{\tilde{\xi}(\theta_0,\theta_1)}
\def\chiitau2{P^i_{\theta}}

\def\myK{K}

\def\tauu{\theta}

\def\theP{P}

\def\thegamma{\gamma}

\def\cA{\mathfrak{A}}\def\cA{\mathcal{A}}
\def\cB{\mathcal{B}}

\def\cD{\mathcal{D}}
\def\cE{\mathcal{E}}
\def\cF{\mathfrak{F}}
\def\cG{\mathfrak{G}}

\def\cL{\mathcal{L}}
\def\cM{\mathcal{M}}

\def\cO{\mathfraka{O}}\def\cO{\mathcal{O}}
\def\cP{\mathfrak{P}}\def\cP{\mathcal{P}}
\def\cQ{\mathcal{Q}}

\def\bX{\mathbb{X}}

\def\emph{}

\def\tb{{\bar{\theta}}}
\def\thisR{\thisR }

\def\textsl{}

\def\proof{\noindent {\it {\textbf{Proof}}}.$\;\,$}
\def\finproof {\hfill $\Box$ \vskip 5 pt }\def\finproof{\rule{4pt}{6pt}}

\def \Int{\displaystyle\int}

\def \be{\begin{eqnarray}}
\def \ee{\end{eqnarray}}
\def \b*{\begin{eqnarray*}}
\def \e*{\end{eqnarray*}}

\def \[{[\,\!\![}
\def \]{]\,\!\!]}

\def \1{{\bf 1}}

\def\F{{\cal T}}\def\F{{\mathfrak{G}}}\def\F{{\mathfrak{F}}}

\def\N{N}\def\N{{\mathbb N}}

\def\I{\mathds{1}}

\def\cL{\mathcal{L}}

\def\hat{\widehat}

\def\mym{m}

\def\xitau2{\xi_{(\tauu)}}\def\xitau2{\xi}
\def\xiitau2{\xi^i_{(\tauu)}}\def\xiitau2{\xi^i}
\def\xintau2{\tilde{\xi}_{(\tauu)}}\def\xintau2{\tilde{\xi}}

\def\Ltau2{\tilde{\xi}(\theta_{\mym},\theta_1)}
\def\chiitau2{P^i_{\theta}}

\def\mynu{\nu}

\def\mylambda{\tilde\lambda}\def\mylambda{\tilde\gamma}

\def\mym{-}

\def\theh{h}

\def\tilde{\widetilde}
\def\cadlag{c\`adl\`ag\xspace}

\def\bX{b}

\def\sigmaX{\sigma}
\def\emph{}

\def\representative{reduction\xspace}

\def\G{{\mathfrak{G}}}\def\G{{\mathfrak{F}}}
\def\ff{\ff}

\def\ff{\ff^{\star}}
\def\F{{\mathfrak{F}}^{\star}}

\def\gg{\gg}\def\gg{{\mathfrak{G}}}
\def\G{{\mathfrak{G}}}
\def\ff{\ff}\def\ff{{\mathfrak{F}}}
\def\F{{\mathfrak{F}}}
\def\hh{\mathfrak{H}}

\def\E{\mathbb{E}}

\def\Rf{\mathfrak{R}}\def\Rf{\Rf}\def\Rf{R_{f}}\def\Rf{\bar{R}_{b}}

\def\thisR{R}\def\thisR{\thisR }\def\thisR{R_b}

\def\tauu{\vartheta}
\def\bthet a{\tau}
\def\tb{{\bar{\theta}}}\def\tb{\tauu}

\def\bthet a{\vartheta}
\def\tb{\tauu}\def\tb{{\bar{\theta}}}

\def\paragraph{\noindent\textbf}

\def\theG{\Gamma}

\def\thec{c}

\def\qqq{\quad\quad\quad}
\def\qq{\quad\quad}

\def\thec{c}
\def\theN{N}

\def\passhortciteN{\citeN}

\def\themu{\mu}

\def\bthet{\boldsymbol\tau}\def\bthet{\mathbf{k}}

\def\thisR{R}

\def\qr#1{\eqref{#1}}
\def\fr#1{Fig.~\ref{#1}}
\def\reft#1{Tab.~\ref{#1}}
\def\sr#1{Sect.~\ref{#1}}
\def\pr#1{Part~\ref{#1}}
\def\refc#1{Chapter~\ref{#1}}

\newcommand{\indi}[1]{\I_{\{{#1}\}}}

\newcommand{\bethe}{\bt}
\newcommand{\ethe}{\et}

\newcommand{\iend}{\end{itemize}}
\newcommand{\ok}{\rule{4pt}{6pt}}\renewcommand{\ok}{\finproof}
\newcommand{\desb}{\begin{description}}
\newcommand{\dese}{\end{description}}

\newcommand{\dcb}{\begin{array}{lll}}
\newcommand{\dce}{\end{array}}
\newcommand{\ebe}{\begin{enumerate}[1)]\setlength{\baselineskip}{13pt}\setlength{\parskip}{5pt}}
\newcommand{\dbe}{\end{enumerate}}

\newcommand{\ibegin}{\begin{itemize}\setlength{\baselineskip}{19pt}\setlength{\parskip}{7pt}}

\makeatletter
\newenvironment{systeme*}{\def\arraystretch{1.2}\arraycolsep=1.4pt\left\{\begin{array}{l}}{\end{array}\right.}

\makeatother

\def\thisG{G}

\def\croc#1{{\langle#1\rangle}}

\def\bal{\begin{aligned}}
\def\eal{\end{aligned}}
\newcommand{\beql}[1]{\beqa\begin{aligned}\label{#1}}
\newcommand{\eeql}{\eal\eeqa}

\newcommand{\bel}{\begin{eqnarray*}\begin{aligned}}
\newcommand{\eel}{\eal\end{eqnarray*}}
 
\def\prad{\otimes}\def\prad{\times}

\newcommand{\eqdef}{\mathrel{\mathop:}=}

\usepackage{hyperref}
\hypersetup{
 colorlinks  = true, %
 urlcolor   = black, %
 linkcolor  = black, %
 citecolor  = black %
}
\makeatletter
\renewcommand\@makefnmark{\hbox{\@textsuperscript{\normalfont\color{blue}\@thefnmark}}}
\renewcommand\@makefntext[1]{%
 \parindent 1em\noindent
      \hb@xt@1.8em{%
        \hss\@textsuperscript{\normalfont\@thefnmark}}#1}
\makeatother
\begin{document}

\def\theW{W}
\def\ftime{\tau}
\def\gtime{\vartheta}
\def\thisg{g}
\def\thisG{G}
\def\wG{\widehat{\mathtt{G}}}\def\wG{\widehat{G}}
\def\thisU{U}
\def\thisW{W}
\def\theG{G}
\def\ttG{\mathtt{G}}
\def\theU{U}\def\theU{\thesigma}\def\theU{\ftime }\def\theU{\theta}
\def\tS{S}
\def\tSigma{\Sigma}\def\tSigma{\ttS}
\def\theM{M}\def\theM{R}
\def\theN{Q}\def\theN{N}
\def\tA{\mathtt{A}}\def\tA{A}
\def\theY{Y}\def\theY{X}
\def\theV{V}\def\theV{\ftime}\def\theV{\nu}\def\theV{\sigma}
\def\theZ{Z}\def\theZ{\Zb}
\def\Zb{\mathsf{Z}}
\def\thatV{V}
\def\then{n}
\def\thisY{Y}
\def\fo{\ff.o}\def\fo{o}
\def\fp{\ff.p}\def\fp{p}
\def\thisM{R}\def\thisM{M}
\def\otherM{\rho}\def\otherM{\RZ}
\def\thisN{N}\def\thisN{P}
\def\thesigma{\theU}
\def\tJ{\mathtt{J}}\def\tJ{J}
\def\tH{\mathtt{H}}\def\tH{H}
\def\tS{\mathtt{S}}\def\tS{S}
\def\ther{r}\def\ther{\varrho}\def\ther{\mathtt{R}}\def\ther{\mathtt{M}}
\def\unrho{\rho}\def\unrho{\mathsf{R}}\def\unrho{R^{\prime}}\def\unrho{M^{\prime}}\def\unrho{M^{\star}}
\def\RX{R}\def\RX{M}\def\RX{M^X}
\def\RY{R^{\prime}}\def\RY{R^{\prime\prime}}\def\RY{M^{\prime\prime}}\def\RY{M^Y}
\def\RZ{R^{\prime\prime}}\def\RZ{\mathsf{R}}\def\RZ{\mathsf{M}}\def\RZ{\mathcal{M}}\def\RZ{M^\star}\def\RZ{M^Z}
\def\monp{\mathsf{p}}\def\monp{p}
\def\monq{\mathsf{q}}\def\monq{q}
\def\thp{\underline{\thisN}}\def\thp{p}
\def\thrh{p}\def\thrh{\rho}
\def\bo{\bar{o}}\def\bo{{\bf o}}
\def\tI{\mathbb{I}}
\def\qr#1{\eqref{#1}}
\def\fr#1{Fig.~\ref{#1}}
\def\reft#1{Tab.~\ref{#1}}
\def\sr#1{Sect.~\ref{#1}}
\def\ar#1{Appendix~\ref{#1}}
\def\pr#1{Part~\ref{#1}}
\def\refc#1{Chapter~\ref{#1}}
\def\Jt{J^\star}\def\Jt{\mathsf{J}^\star}
\def\gr{[\![}\def\gr{[}
\def\rg{]\!]}\def\rg{]}
\def\dHt{\boldsymbol\delta_\gtime(dt)}\def\dHt{dH_t}

\def\thatF{F}\def\thatF{\mathtt{G}}

\def\thisGp{\theF}
\def\thisGm{\overline{\theF}}\def\thisGm{\theF^{\prime}}
\def\theF{\mathtt{F}}
\def\subset{\subseteq}

\def\thealpha{\alpha}
\def\cM{Q}\def\cM{\mathtt{M}}\def\cM{\mathtt{Q}}
\def\Alpha{\mathtt{D}}
\def\pAlph{\mathtt{A}}
\def\ttS{\mathtt{S}}
\def\pSigma{\hat{\tSigma}}\def\pSigma{{^{p}}\Sigma}\def\pSigma{{^{p}}\!\ttS}
\def\PX{\Psi}
\def\PZ{\Psi'}
\def\PV{\Psi''}
\def\HV{\Phi^{\thatV}}\def\HV{\Phi''}
\def\PX{\Psi^X}
\def\PZ{\Psi^Z}
\def\PV{\Psi^V}
\def\HV{\Phi^V}
\def\HX{\Phi^X}
\def\HZ{\Phi^Z}
\def\HP{\Phi^P}\def\HP{}
\def\HPi{(\Phi^P)^{-1}}\def\HPi{}
\def\MX{N^X}
\def\MZ{N^Z}\def\MZ{M^Z}
\def\MV{P^V}
\def\RX{M^X}
\def\RZ{M^Z}

\def\Es{\mathcal{E}^{\star}}

\def\monqs{\mathcal{E}(\ind_{\{\pSigma >0\}}\frac{1}{\pSigma }\centerdot \cM)}\def\monqs{\mathcal{E}(\frac{1}{\pSigma }\centerdot \cM)}
\def\In{\cup_n[0,\varsigma_n]}\def\In{\{\ttS_{-}>0\}}
\def\Int{\cup_n[0,T\wedge\varsigma_n]}\def\Int{\In}
\def\pIn{(\In)}
\def\pInt{(\Int)}
\def\therho{\rho}\def\therho{\rho}\def\therho{\mathfrak{q}}\def\therho{\mathsf{q}}

\def\ttJ{{\ind_{[0,\tau)}}}\def\ttJ{\mathtt{J}}
\def\DE{(\Delta^{e} )_{e\in E}}\def\DE{\Delta}
\def\DP{(\Delta^{e'} )_{e'\in E}}\def\DP{\Delta}
\def\DT{(\widetilde{\Delta}^e)_{e\in E}}\def\DT{\widetilde{\Delta}}

\def\fact{\mathbb{Q}(d\omega)\prad}\def\fact{q^e_s\gamma_s}

\def\thec{c}
\def\thegamma{\gamma}
\def\thejX{\thej}
\def\sigmaX{\sigma} 
\def\bX{\beta}\def\bX{b}
\def\partialx{\partial}
\def\partialxx{\partial^2\!}
\def\thej{j}\def\thej{\psi}
\def\Pssi{\Psi}
\def\xd{x}

\def\x{x}
\def\e{e}
\def\de{de} 
\def\Xb{X}
\def\thatj{j}

\def\sr#1{Section~\ref{#1}}

\def\ttH{\mathtt{H}}\def\ttH{{\ind_{[\tau,\infty)}}}
\def\Jt{\ind_{(0,\tau]}}
\def\theh{g}
\def\theP{M'}\def\theP{P}
 
 \def\MRP{martingale representation property\xspace}
 \def\themu{\mu}\def\themu{\eta}\def\themu{\beta}
 \def\mylambda{\gamma}
 
\def\tb{\bar{\tau}}

\def\theV{\theta'}
\def\myK{L}\def\myK{M'}

\def\mathtt{\mathsf}
\def\sta{\star}\def\sta{*}

\def\g{\textcolor{ForestGreen}}
\def\r{}\def\r{\textcolor{red}}
\def\b{\textcolor{blue}}\def\b{}

\author{St\'ephane Cr\'epey\footnote{{%
\texttt{stephane.crepey@lpsm.paris}. This research has benefited from the support of the Chair \textit{Capital Markets Tomorrow: Modeling and Computational Issues  under the aegis of the Institut Europlace de Finance},  a joint initiative of Laboratoire de Probabilit\'es, Statistique et Mod\'elisation (LPSM) / Universit\'e Paris Cit\'e and  Cr\'edit Agricole CIB. The author is grateful to Shiqi Song for his contributions to a preliminary version of this work and to Martin Schweizer and Monique Jeanblanc for precious comments.}}
\\Laboratoire de Probabilit\'es, Statistique et Mod\'elisation (LPSM), \\Sorbonne Universit\'e et Universit\'e
Paris Cit\'e, CNRS UMR 8001}

\def\gammap{\gamma'}\def\gammap{\gamma}
\def\varphi{\gammap }
\def\Et{\tilde{\E}}\def\Et{\E^{\P}}\def\Et{\E'}
\def\thevarrho{\varrho}\def\thevarrho{\eta}
\def\cE{\mathcal{E}}
\title{Invariance Times Transfer Properties}

\maketitle

\begin{abstract} 
Invariance times are stopping times $\tau$ such that local martingales with respect to some reduced filtration and an equivalently changed probability measure, stopped before $\tau$, are local martingales with respect to the original model filtration and probability measure. They arise naturally for modeling the default time of a dealer bank, in
the mathematical finance context of counterparty credit risk. Assuming an invariance time endowed with an intensity and a 
positive Az\'ema supermartingale, 
this work establishes a dictionary relating the semimartingale calculi in the 
original and reduced stochastic bases, regarding conditional expectations,
martingales, 
stochastic integrals, 
random measure stochastic integrals,
martingale representation properties,
semimartingale 
characteristics, 
Markov properties, transition semigroups and infinitesimal generators,
and solutions of backward stochastic differential equations.
\end{abstract}
\begin{keywords}
Progressive enlargement of filtration, invariance time, semimartingale calculus, 
Markov process, backward stochastic differential equation,
counterparty risk, credit risk.
\end{keywords}

\vspace{2mm}
\noindent
\textbf{Mathematics Subject Classification:} 60G07, 60G44, 60H10, 	91G20, 91G40.

\section{Introduction}\label{s:intro}

This paper is about a concept in progressive enlargement of filtrations called invariance times introduced by \citet{CrepeySong15c}. Progressive enlargement of filtrations refers to a situation where two filtrations are involved, a smaller and a bigger one, the bigger one making a certain random time in the small one, $\tau$, a stopping time. $\tau$ is called an invariance time
  when local martingales $X$ in the small filtration, once stopped before $\tau$, so
 $X^{\tau -}\eqdef  X\ind_{[0,\tau )} +X_{\tau -} \ind_{[\tau ,+\infty)},$
are local martingales in the large one. However, this is not necessarily required to hold under the original probability measure, but only under a possibly modified one. In other words, there exists a measure change that ``compensates'' the change of filtration. The basic situation, called immersion, is when local martingales in the small filtration do not jump at $\tau$ and are local martingales in the large one, for the original probability measure (so no measure change is required). But 
  immersion is tantamount to a certain form of independence between $\tau$ and the small filtration \citep*[Lemma 3.2.1(ii)-(iii)]{Bielecki2009}. This is too restrictive for applications to credit risk in finance, which is a lot about the dependence between the market risk represented by the small filtration and the default time $\tau$ (especially adverse dependence, dubbed ``wrong way risk''). \citep{CrepeySong15FS,CrepeySong16b} show that invariance times offer a much more flexible framework in this regard. 
  
  The ``stopping before $\tau$'' feature in the above resonates particularly well with a particular application, namely the pricing of the implications for a bank of its own default time $\tau$. Indeed, for the shareholders of the bank, only the pre-default cash flows matter, hence the corresponding pricing equations are stopped before $\tau$. These equations are known as the XVA equations, where X is a catch-up letter to be replaced by C for credit, F for funding or K for capital, while VA stands for valuation adjustment \citep[Eqns. (2.12), (2.13), and (2.17)]{Crepey21}. The setup of invariance times arises as one which is together flexible  enough in terms of credit-market dependence, whilst being amenable to an elegant solution of the XVA equations, by reduction to simpler equations stated with respect to a smaller filtration in which the default risk of the bank only appears through its intensity. However, to perform this reduction rigorously, one needs to relate the stochastic calculi in the small and the large filtrations, regarding conditional expectations,
martingales, 
stochastic integrals, 
random measure stochastic integrals,
martingale representation properties,
semimartingale 
characteristics, 
Markov properties, transition semigroups and infinitesimal generators,
and (eventually) solutions of backward stochastic differential equations (BSDEs).
The elaboration of the corresponding transfer properties is the contribution of this paper. 
Section \ref{ss:C} sets the stage.
 The conditional expectations transfer formulas of Section \ref{ss:etf} underlie most of the subsequent developments
 of Sections \ref{ss:ex-suite}--\ref{s:markov}
  that, apart from this common base, are quite independent from each other. Hence the reader can cherrypick freely among these. The 
  BSDE Section \ref{s:bsde} puts 
more or less 
  everything together.
 Section \ref{s:concl} concludes.
Our conditional expectations transfer formulas present a resemblance with earlier formulas stated in terms of a singular measure change but no change of filtration \citep*[Theorem 1]{CollinDufresneGoldsteinHugonnier2004}.  The connection between the two approaches is illustrated in the end of the paper.

\subsection{Standing Notation and Terminology}\label{ss:san}

The real line and half-line are denoted by $\R$ and $\R_+$;
$ | \cdot|$ denotes any Euclidean norm (in the dimension of its argument), $\cdot^\top$ means vector transposition;
$\cB (E)$ denotes the Borel $\sigma$ algebra on a metrizable space $E$;
$\boldsymbol\lambda$ is the Lebesgue measure on $\R_+$,
$\boldsymbol\delta_a$ represents
a Dirac measure at a point $a$.

Unless otherwise stated, a function (or process) is real valued; order relationships between random variables (respectively processes) are meant almost surely (respectively in the indistinguishable sense); a time interval is random  (in particular, the graph of a random time $\theta$ is simply written $\gr \theta \rg$\index{Ã©@$\gr \theta\rg$}).
We do not explicitly mention the domain of definition of a function 
when it is implied by the measurability, e.g.~we write
``a $\cB ({\R})$ measurable function $h$ (or $h(x)$)'' rather than
``a $\cB ({\R})$ measurable function $h$ defined on ${\R}$''.
For a function $h(\omega,x)$ defined on a product space $\Omega\times E$, we write $h(x)$ (or $h_t$ in the case of a stochastic process), without $\omega$.

We use the terminology of the general theory of processes and of filtrations as given in the books by \citeN{DM75}
and \passhortciteN{HeWangYan92}. 
For any semimartingale, always taken in a \cadlag version in this work, $X$, and for any predictable $X$ integrable process $L$, the corresponding stochastic integral is denoted by $\int_0^{\cdot} L_t dX_t=\int_{(0,\cdot]}L_t dX_t=L\centerdot X$, with the precedence convention \index{\'e@$\centerdot$}$KL\centerdot X=(KL)\centerdot X$ if $K$ is another predictable process such that $KL$ is $X$ integrable. The stochastic exponential of a semimartingale $X$ is denoted by $\cE(X)$. By drift of a special semimartingale (i.e.~a semimartingale with locally integrable jumps), 
we mean the finite variation predictable part of its canonical Doob-Meyer decomposition. Stochastic integrals of random functions with respect to jump measures and their compensations are meant in the sense of \citeN{Jacod1979}, to which
we also borrow the usage of including the optionality with respect to a reference filtration in the definition of an integer valued random measure. 
Random measure stochastic integrals and transform of measures by densities are respectively denoted by ``$*$'' and ``${\cdot}$''.  We denote by
$\cP (\hh )$\index{p@$\cP (\hh )$} and 
$\cO(\hh )$\index{o@$\cO(\hh )$} the predictable and optional
$\sigma$ fields with respect to a filtration $\hh $. 

For any random time $\theU $ and c\`adl\`ag process $X$, $\Delta_\theU X$ represents the jump of $X$ at $\theU .$
We use the convention that \index{\'e@${\cdot}_{0-}$}$X_{0-}=X_0$ (hence $\Delta_0 X=0$) 
and we write \index{\'e@${\cdot}^{\theU -}$}$X^{\theU}$ and \index{\'e@${\cdot}^{\theU -}$}$X^{\theU -}$ for the processes $X$ stopped at $\theU$ and before $\theU $, i.e.
\beql{right}
X^{\theU }=X\ind_{[0,\theU )} +X_{\theU} \ind_{[\theU ,+\infty)} \sp X^{\theU -}=X\ind_{[0,\theU )} +X_{\theU -} \ind_{[\theU ,+\infty)}.
\eeql
The process $X$ is said to be stopped at $\theU$, respectively before $\theU$, if
$X=X^{\theU}$, respectively $X=X^{\theU-}$.
We call compensator of a stopping time $\theta$ the compensator of $\ind_{[\theta,\infty)}$.
We say that $\theta$ has an intensity $\gamma$ if
 $\theta$ is (strictly) positive and that its
compensator is given as 
$\gamma\centerdot\boldsymbol\lambda$, for some predictable process $\gamma$ (vanishing beyond time $\theta$).
For any event $A$, we denote by $
\theta _A$ the stopping time $\ind_A
\theta + \ind_{A^c} \infty.$

\section{Setup}\label{ss:C}

We work on a space $\Omega$ equipped with a $\sigma$ field \index{a@$\mathfrak{A}$}$\cA$, a probability measure \index{q@$\mathbb{Q}$}$\mathbb{Q}$ on $\cA$, and a filtration $\gg=(\mathfrak{G}_t)_{t\in\mathbb{R}_+}$
\index{g@$\gg$} of sub-$\sigma$ fields of $\cA$ satisfying the usual conditions.
We are given a 
positive (nonnecessarily finite) 
$\gg$ stopping time \index{t@$\ftime$}$\tau$ and 
a subfiltration 
\index{f@$\ff$} 
$\ff=(\mathfrak{F}_t)_{t\in\mathbb{R}_+}$ 
of $\gg$ 
satisfying the usual conditions, with $\ff$ optional and predictable projections denoted by \index{o@${^{\fo }}$}${^{\fo }}\!\cdot$ and \index{p@${^{p}}$}${^{p}}\!\cdot$.
We consider the 
progressive enlargement of filtrations setup\footnote{recalling from \citeN[Section 2.1]{CrepeySong15c}
that the proofs of the classical progressive enlargement of filtration results in \citeN{JeulinYor78} or Chapitre 20 in \citeN{DellacherieMaisonneuveMeyer92}, although only stated and proved in the specific setup of Example \ref{ex:classicalpef},  
all work in the extended setup \eqref{e:test}.} defined by the condition
that
\begin{equation}\label{e:test}
\forall t\ge 0 \mbox{ and } B \in {\mathfrak{G}_t }\sp \exists B' \in\mathfrak{F}_t\mbox{ such that } B \cap\{t<\ftime \}=B' \cap\{t<\ftime \}.\end{equation}
\bex\label{ex:classicalpef}
This holds in particular (but not only, see \sr{ss:ex}) in the classical progressive enlargement of filtration setup
$$\G_t= \F_{t} \vee    \sigma(\tau
    \wedge t)\vee \sigma(\{\tau
    > t\}) \sp t\in\R_+,$$
i.e.~when $\gg$ is the smallest filtration larger than $\ff$ making $\tau$ a stopping time.~\finproof
 \eex
Equivalently to \eqref{e:test}\footnote{see \citeN[Eqn. (2.1)]{CrepeySong15c}.},
any ${\gg}$ predictable (resp.~optional) process $L$ admits an ${\ff}$ predictable (resp.~optional) process $L'$\footnote{the notation $\cdot'$  common to the predictable and optional projections is typically not an issue in practice (whenever useful we write explicitly ``predictable'' or ``optional'').}, dubbed $\ff$ predictable (resp.~optional)  \representative 
of $L$, such that $\ind_{(0,\tau]} L = \ind_{(0,\tau]} L' $ (resp.~$\ind_{[0,\tau)}L=\ind_{[0,\tau)} L')$.  In particular, for any $\gg$ stopping time $\theU$, there exists an $\ff$ stopping time $\theV$, dubbed $\ff$ \representative of $\theU$, such that ${\{\theU<\ftime \}}={\{\theV<\ftime \}}\subseteq {\{\theU=\theV\}}$.

Given a positive constant $T$, 
we work henceforth under the following condition, introduced with its first consequences in \citeN[Sections 4--6]{CrepeyElie16}, and which is explored systematically in this work.\\

\noindent\textbf{Condition \index{a@(C)}(C)}. $\tau$ has a $(\gg,\mathbb{Q})$ intensity, the 
Az\'ema supermartingale $\ttS ={^{\fo }}\! (\ind_{[0,\tau)})$ of $\tau$
satisfies $\mathtt{S}_T >0$ almost surely, and 
\beql{e:A}& 
\mbox{there exists
a 
probability measure $\mathbb{P}$
equivalent to $\mathbb{Q}$ on $\cF_T$, called invariance}\\&\mbox{probability measure,
such that, for any $({\ff},\mathbb{P})$ local martingale $\thisN$,
}\\&\mbox{$\thisN^{\tau-}$
is a $(\gg,\mathbb{Q})$ local martingale on $[0,T]$.}~\finproof \eeql

\noindent
Unless explicitly stated, the reference probability measure is the original measure $\Q$.

The conjunction of \eqref{e:test} and \eqref{e:A} corresponds to the notion of invariance time $\tau$.
Hence all the results of \citeN{CrepeySong15c} are applicable in this work, 
sometimes in a stronger version due to the additional assumptions embodied in the first line of the condition (C). In particular:
\bl \label{lem:pre} 
Under the condition (C),
{\rm
\hfill\break
 \textbf{(i)}} 
$\{ \ttS_->0\}=\{\pSigma >0\}=\{ \ttS  >0\}\supseteq [0,T];$
{\rm\hfill\break \textbf{(ii)}} two 
$\ff$ optional processes that coincide before $\tau$ coincide on $[0,T]$;  in particular,
predictable and optional reductions are uniquely defined on $[0,T]$;
{\rm\hfill\break \textbf{(iii)}}  
invariance probability measures $\P$ are uniquely determined on $\cF_T$, with  
$(\ff,\Q)$ density process  
\beql{e:conda}\cQ:=\monqs\mbox{  on $ [0,T]$},\eeql
a positive $\ff$ martingale on $ [0,T].$ 
\el
\proof {\rm \textbf{(i)}} 
By Theorem 3.7 in \citeN{CrepeySong15c}, in the case of an invariance time $\tau$ endowed with a $(\gg,\Q)$ intensity,
$\{ \ttS_->0\}=\{\pSigma >0\}=\{ \ttS  >0\}.$
Under the additional
assumption ``$\mathtt{S}_T >0$ a.s.'' that is postulated in the condition (C), we can add 
`` $\supseteq [0,T]$'' .
{\rm\hfill\break \textbf{(ii)}} The first part is Lemma 2.3 in \citeN{CrepeySong15c}, which readily implies the second part.
 {\rm\hfill\break \textbf{(iii)}} By Theorem 3.2 in \citeN{CrepeySong15c} and (i).~\finproof\\

\noindent
As $\P$ is only used for computations in $\ff$ on $[0,T]$, it only matters on $\ff_T$.
Hence, in view of Lemma \ref{lem:pre}(iii), 
we can talk of ``the invariance probability measure $\P$'' in our setup. Moreover,
by reduction, we  may and do assume that the $(\gg,\mathbb{Q})$  intensity of $\tau$ is of the form {$\gamma\ind_{(0,\tau]}$}, for an 
$\ff$ predictable process $\gamma$ uniquely defined on $[0,T]$,   by Lemma \ref{lem:pre}(ii).
We write $\Gamma=\int_0^{\cdot}\gamma_s ds, $ so that 
$\Gamma^{\tau}$ is the $(\gg,\mathbb{Q})$  compensator of $\tau$.

\section{Conditional Expectation Transfer Formulas}\label{ss:etf}

The $(\G_t,\Q)$ and $(\F_t,\P)$  conditional expectations are denoted by
$\E_t$ and $\Et_t$ and we drop the index $t$ at time 0.

The following result, the unconditional version of which corresponds to Theorem 4.1 in \citeN{CrepeyElie16}, 
provides an extension of classical results
(see e.g.~\citeN[Chapter 3]{Bielecki2009}) beyond
the basic immersion setup where $(\ff,\P=\Q)$ local martingales are $(\gg,\Q)$ local martingales without jump at $\tau$.

\bt
\label{t:exptransf}
For any constant $t\in [0,T]$, any
$[t,T]$ valued 
$\ff$ stopping time $\sigma$, any $\mathfrak{F}_\sigma$ measurable nonnegative random variable $\chi$, any $\ff$ predictable nonnegative process $K,$ 
and any $\ff$ optional nondecreasing process $A$ starting from 0,
we have, on $\{ t< \tau\}$,
\begin{eqnarray}\label{t:etf}
&&\mathbb{E}_t [\chi\ind_{\{\sigma<\tau\}} ]
= 
 \Et_t[\chi e^{- (\Gamma_\sigma-\Gamma_t) } ]
\sp
\\&&\label{t:etfp}
\mathbb{E}_t [K_{\tau}\ind_{\{ \tau\leq T\}} ]
=
\Et_t \Big[\int_t^T \ K_{s}e^{- (\Gamma_s -\Gamma_t)}\varphi_s \ ds  \Big] , 
\\ \label{t:etfpp}
&&
\mathbb{E}_t [A^{\tau-}_T- A^{\tau-}_t  ]
=
\ind_{\{ t< \tau\}}\Et_t \Big[\int_t^T \ e^{-(\Gamma_s-\Gamma_t)} \ dA_{s}\Big].
\end{eqnarray} 
\et
\proof 
Consider the $\ff$ canonical Doob-Meyer decomposition $\ttS=\cM-\Alpha$ of  $\ttS$,
where $\cM$ {(with $\cM_0=\ttS_0=1$)} and
$\Alpha$  {(with $\Alpha_0=0$)} are
the $\ff$ local martingale component
and  the $\ff$
drift of $\ttS.$  
By Lemma 2.2 5) in \citeN{CrepeySong15c} and Lemma \ref{lem:pre}(i), in the present setup where $\tau$ is positive, so that $\ttS_0=1$, and Lemma \ref{lem:pre}(i) is satisfied, $\ttS$ admits the multiplicative decomposition
\beql{e:multdec} \ttS=  \mathcal{Q}\mathcal{D}
\mbox{ on } [0,T],\eeql
where $\mathcal{Q}$ is the $\ff$ martingale \eqref{e:conda} on $[0,T]$
and
$\mathcal{D}=\cE(-\frac{1}{\ttS_-}\centerdot \Alpha)$
is an $\ff$ predictable nonincreasing process on $[0,T]$. 

For any $B\in \G_t$ and $B'$ associated with $B$ as in \qr{e:test},
we then have by definition of $\mathtt{S}_s=\Q(\tau
>s\,|\,\cF_s), s\ge 0$,
 and $\mathfrak{F}_\sigma$ measurability of $\chi$ 
 (using also the tower rule and recalling the assumption ${\ttS}_T>0$ in the condition (C)):
\begin{eqnarray}\nonumber\begin{aligned}
& \mathbb{E}\left[\ind_{\{ t< \tau\}} \mathbb{E}\big( \chi {\ttS}_\sigma / {\ttS}_t\big|\cF_t\big)   \ind_{B} \right]
= \mathbb{E}\left[{\ttS}_t \mathbb{E}\big(\ind_{B'} \chi {\ttS}_\sigma / {\ttS}_t  |\cF_t \big)   \right] 
=\mathbb{E}\left[ \chi {\ttS}_\sigma  \ind_{B'}   \right]
=\mathbb{E}\left[ \chi \ind_{\{ \sigma< \tau\}} \ind_B  \right] .
\end{aligned}\end{eqnarray} 
Hence
\beql{e:tempo}
& \ind_{\{ t< \tau\}} \mathbb{E}\big( \chi {\ttS}_\sigma / {\ttS}_t\big|\cF_t\big)    
=   \mathbb{E}\big(\ind_{\{  \sigma< \tau\}}   \chi   |\cG_t \big) .
\eeql 
Then 
\qr{e:multdec}, 
under the assumption ${\ttS}_T>0$ a.s., yields 
\beql{e:fpgqprel}
&\mathbb{E}\big( \chi {\ttS}_\sigma / {\ttS}_t\big|\cF_t\big)    =
\mathbb{E}\Big( \chi \cQ_\sigma \cD_\sigma /  \big(\cQ_t \cD_t\big)\Big|\cF_t\Big)  =\Et\Big[\chi    \cD_\sigma / \cD_t \Big|\cF_t \Big], 
\eeql
by Lemma \ref{lem:pre}(iii) and the conditional Bayes formula corresponding to the $\Q$-to-$\P$ density process $\cQ$ on $[0,T]$. 
Moreover,
by \citeN[Lemma A.1]{CrepeySong15c} and Lemma \ref{lem:pre}(i),
$\mathtt{D}$ is continuous and
\beql{e:inta}\cD^{\pm 1}=\cE(\pm\frac{1}{\ttS_-}{\centerdot}\mathtt{D})=e^{\pm\frac{1}{\ttS_-}{\centerdot}\mathtt{D}}\sp \frac{1}{\ttS_-}{\centerdot}\mathtt{D}=\gamma'  {\centerdot}\boldsymbol{\lambda}\eeql
hold
on $ [0,T]$.
Recalling that $\Gamma=\gamma  {\centerdot}\boldsymbol{\lambda}$,
\eqref{e:tempo}--\eqref{e:inta} yield \qr{t:etf}. 

For \qr{t:etfp}, we compute, on $ \{ t< \tau\} $,
\bel
&
\mathbb{E}_t[K_{\tau}\ind_{\{\tau\leq T\}}]
=
\mathbb{E}_t \Big[\int_t^T K_{s}\ind_{\{s\leq \tau\}}\varphi_s\ ds\Big]
=
\int_t^T  \mathbb{E}_t[K_{s}\ind_{\{s< \tau\}}\varphi_s ]\ ds\\
&\qqq= 
\int_t^T \ \Et_t[K_{s}e^{-(\Gamma_s-\Gamma_t)}\varphi_s  ]\ ds
=
\Et_t \Big[\int_t^T \ K_{s}e^{-(\Gamma_s-\Gamma_t)}\varphi_s \ ds\Big],
\eel
where \qr{t:etf} was used for passing to the second line.

Regarding \qr{t:etfpp}, an application of \qr{t:etfp} yields  (still on $ \{ t< \tau\} $)
\bel
&
\mathbb{E}_t[(A_{\tau-}-A_t)\ind_{\{\tau\leq T\}}]
=
\Et_t \Big[\int_t^T \ (A_{s}-A_t)e^{-(\Gamma_s-\Gamma_t)}\varphi_s \ ds\Big]
\\&\qqq
=
-\Et_t[(A_{T}-A_t)e^{-(\Gamma_T -\Gamma_t)}]
+\Et_t \Big[\int_t^T \ e^{-(\Gamma_s-\Gamma_t)} \ dA_{s}\Big].
\eel
Using \qr{t:etf}, we deduce
\bel
\mathbb{E}_t[(A^{\tau-}_T&-A_t)]
= 
\mathbb{E}_t[(A_{T}-A_t)\ind_{\{T<\tau\}}]+\mathbb{E}_t[(A_{\tau-}-A_t)\ind_{\{\tau\leq T\}}]\\
&= 
\Et_t[(A_T-A_t) e^{-(\Gamma_T -\Gamma_t)}]
-\Et_t[(A_{T}-A_t)e^{-(\Gamma_T -\Gamma_t)}]
+\Et_t \Big[\int_t^T \ e^{-(\Gamma_s-\Gamma_t)} \ dA_{s}\Big]\\
&= 
\Et_t \Big[\int_t^T \ e^{-(\Gamma_s-\Gamma_t)} \ dA_{s}\Big].~\finproof 
\eel
\noindent
See \sr{ss:ex} for a discussion of two
alternatives to the formula \qr{t:etfp} that are known from the mathematical finance literature.

\section{Martingale Transfer Formulas}\label{ss:ex-suite}

We denote by:
\begin{itemize}
\item $\mathcal{M}_{T}(\ff,\mathbb{P})$, the set
 of $(\ff,\mathbb{P})$ local martingales stopped at $T$;
\item $\mathcal{M}_{\tau- \wedge T}(\gg,\mathbb{Q})$, the set
 of $(\gg,\mathbb{Q})$ local martingales stopped at $\tau- \wedge T$, i.e.~before $\tau$
and at $T$;
\item  $\mathcal{M}^c_{T}(\ff,\mathbb{P})$ and $\mathcal{M}^d_{T}(\ff,\mathbb{P})$, respectively $\mathcal{M}^c_{\tau \wedge T}(\gg,\mathbb{Q}),$ and $\mathcal{M}^d_{\tau- \wedge T}(\gg,\mathbb{Q})$, their respective subsets of continuous local martingales and purely discontinuous local martingales.
\end{itemize}
 
\bl\label{l:MN} For any $M,N \in\mathcal{M}_{\tau- \wedge T}(\gg,\mathbb{Q}) $ , we have
\beql{e:MN}
& [M,N]'=[M',N']\mbox{ on } [0,T],
\eeql 
where the quadratic variations $[M,N]$ (with $\ff$ optional reduction $[M,N]'$) and $[M',N']$ are respectively meant in $(\gg,\mathbb{Q})$ and $(\ff,\mathbb{P})$.
\el
\proof
As $M$ and $N$ are stopped before $\tau$,
\beql{e:MNbis}
& [M,N]=[M ,N] ^{\tau-}=[M' ,N'] ^{\tau-}, 
\eeql 
where the quadratic variation $[M',N']$ is meant at this stage in $(\gg,\mathbb{Q})$. But since $(M')^T$ and $(N')^T$ are $\ff$ adapted and that $\P$ and $\Q$ are equivalent on $\cF_T$, on $[0,T]$ this quadratic variation $[M' ,N']$ in $(\gg,\mathbb{Q})$  is the same as the quadratic variation $[M' ,N']$ meant in $(\ff,\mathbb{P})$.
Hence on $[0,T]$ this quadratic variation and $[M,N]'$ coincide before $\tau$, so that they coincide on $[0,T]$, by Lemma \ref{lem:pre}(ii).~\finproof

\bethe\label{invarTH}  
The following bijections hold:
\beql{e:bijdtc}
&\mathcal{M}_T (\ff,\mathbb{P})
\left.
\def\arraystretch{0.5}
\begin{array}{c}\stackrel{\cdot^{\tau-}}{\longrightarrow}\\
\underset{(\cdot^{\prime})^T}{\longleftarrow}  
\end{array} 
\right.    
\mathcal{M}_{\tau- \wedge T}(\gg,\mathbb{Q}),
\\& \mathcal{M}^{c}_T (\ff,\mathbb{P})
\left.
\def\arraystretch{0.5}
\begin{array}{c}\stackrel{\cdot^{\tau}}{\longrightarrow}\\
\underset{(\cdot^{\prime})^T}{\longleftarrow}  
\end{array} 
\right.  
 \mathcal{M}^{c}_{\tau \wedge T}(\gg,\mathbb{Q}),
\\ &\mathcal{M}^{d}_T (\ff,\mathbb{P})
\left.
\def\arraystretch{0.5}
\begin{array}{c}\stackrel{\cdot^{  \tau-}}{\longrightarrow}\\
\underset{(\cdot^{\prime})^T}{\longleftarrow}  
\end{array} \right.  
\mathcal{M}^{d}_{\tau- \wedge T}(\gg,\mathbb{Q}),
\eeql
where $(\cdot^{\prime})^T$ denotes the $\ff$ optional reduction operator composed with stopping at $T$.
\ethe

\proof 
By the converse part in \citeN[Lemma 2.2 4)]{CrepeySong15c} combined with Lemma \ref{lem:pre}(i)-(ii),   
for any $M\in\mathcal{M}_{\tau- \wedge T}(\gg,\mathbb{Q}),$ the process
$\mathtt{S}_-{\centerdot } {M'}+[\mathtt{S},{M'}]$ is an
$(\ff,\mathbb{Q})$ local martingale on $[0,T]$. This in turn implies that
$(M')^T \in \mathcal{M}_T(\ff,\mathbb{P}),$
by Theorem 3.7 in \citeN{CrepeySong15c} and Lemma \ref{lem:pre}(i).
Hence, on $\mathcal{M}_{\tau- \wedge T}(\gg,\mathbb{Q})$, the operator $(\cdot^{\prime})^T$ takes its values  in the space $\mathcal{M}_{T}(\ff,\mathbb{P})$.
Conversely, for any $P\in\mathcal{M}_{T}(\ff,\mathbb{P})$,  the condition (C) yields $P^{\tau-}\in \mathcal{M}_{\tau- \wedge T}(\gg,\mathbb{Q}),$ i.e.
on $\mathcal{M}_T (\ff,\mathbb{P})$ the map $\cdot^{\tau-}$
takes its values in the space
$\mathcal{M}_{\tau- \wedge T}(\gg,\mathbb{Q})$.

To establish the first bijection in \qr{e:bijdtc} it remains to show that $((M^\prime)^T)^{\tau-}=M$ and  
$((P^{\tau-})^\prime)^T=P$ in the above.
As $M$ is stopped before $\tau$ and at $T$, the first identity is trivially true. Regarding the second one,
$(P^{\tau-})'=P$ holds
before $\tau,$ hence on 
$[0,T]  $, by 
Lemma \ref{lem:pre}(ii). Hence $((P^{\tau-})^\prime)^T=P$ holds on $\R_+$ (as $P$ is stopped at $T$).

The second bijection in \qr{e:bijdtc} follows by the same steps, noting that
the reduction of a continuous process $X$ is continuous on $[0,T]$, by Lemma \ref{lem:pre}(ii) applied to the jump process of $X$.

To prove the third bijection, following \shortciteN[Theorem 7.34]{HeWangYan92}, assuming $M\in \mathcal{M}_{\tau- \wedge T}(\gg,\mathbb{Q}),$ we take a $(\gg, \mathbb{Q})$ continuous local martingale $X$ and we consider the bracket $[M, X]$. 
By Lemma \ref{l:MN},
on $[0,T]$, $[M', X']$ is the $\ff$ optional reduction of $[M, X]$. Consequently, using also Lemma \ref{lem:pre}(ii), $[M, X]^{\tau-}=0$ on $[0,T]$ if and only if $[M', X']=0$ on $[0,T]$. The lemma then follows from the first and second bijections in \qr{e:bijdtc}.~\ok

\section{Transfer of Stochastic Integrals in the Sense of Local Martingales}

\bl\label{Tprime}
Let $(\theta_n)_{n\geq 0}$ be a nondecreasing sequence of $\gg$ stopping times tending to infinity.
There exists a nondecreasing sequence $(\sigma_n)_{n\geq 0}$ of $\ff$ stopping times such that $\sigma_n$ tends to infinity and $$
\theta_{n}\wedge T\wedge \tau=\sigma_n\wedge T\wedge \tau.
$$ 
\el

\proof
We compute, using \qr{t:etf}  at $t=0$ for passing to the second line,
\bel
&\Et[\ind_{\{\theta'_n< T\}}e^{-\Gamma_{T}}]
\leq 
\Et[\ind_{\{\theta'_n< T\}}e^{-\Gamma_{\theta'_n}}] 
\\&\qqq=
\mathbb{E}[\ind_{\{\theta'_n< T\}}\ind_{\{\theta'_n< \tau\}}]
=
\mathbb{E}[\ind_{\{\theta_{n}< T\}}\ind_{\{\theta_{n}< \tau\}}]
\rightarrow 0\mbox{ as } n\to \infty. 
\eel
This implies that $\mathbb{P}[\theta'_n<T]\rightarrow 0$. Hence $\mathbb{Q}[\theta'_n<T]\rightarrow 0$, as $\mathbb{P}$ is equivalent to $\mathbb{Q}$ on $\mathfrak{F}_{T}$. The sequence $\sigma_n =(\theta'_n)_{\{\theta'_n<T\}}, n\geq 0,$ satisfies all the desired properties.  \ok\\

\bl\label{integrabilityT}
Let $A$ be a $\gg$ adapted nondecreasing \cadlag process. The process $A^{\tau-}$ is $(\gg,\mathbb{Q})$ locally integrable on $[0,T]$ if and only if $A'$ is $(\ff,\mathbb{P})$ locally integrable on $[0,T]$.
\el

\proof First note from 
\citep[Lemma 6.10]{Song16a} (with $\mathtt{S}_T>0$ a.s.~under the condition (C))
that $A'$ is a nondecreasing process on $[0,T]$.
Let $(\theta_{n})_{n\geq 0}$ be a nondecreasing sequence of $\gg$ stopping times tending to  infinity. Let $(\sigma_{n})_{n\geq 0}$ be associated with $(\theta_{n})_{n\geq 0}$ as in Lemma \ref{Tprime}.  
We compute
\bel
&\mathbb{E}[\int_0^{\theta_{n}\wedge T}\ind_{\{s<\tau\}} e^{\Gamma_s} dA_s\ ]
=
\mathbb{E}[\int_0^{\theta_{n}\wedge T\wedge \tau}\ind_{\{s<\tau\}}e^{\Gamma_s}dA_s\ ]\\&\qqq
=
\mathbb{E}[\int_0^{\sigma_{n}\wedge T\wedge \tau}\ind_{\{s<\tau\}}e^{\Gamma_s}dA'_s\ ]
=
\mathbb{E}[\int_0^{\sigma_{n}\wedge T}\ind_{\{s<\tau\}}e^{\Gamma_s}dA'_s\ ]
=
\Et[A'_{\sigma_{n}\wedge T} ],
\eel
by \qr{t:etfpp} (used at $t=0$).
\b{As $\Gamma$ is continuous, the factor $e^{\Gamma_s} $ can be handled by another layer of localization.}
This implies the result.~\ok

\bethe\label{itf-usualmartingale}
Let 
$W$ be a $(\gg,\mathbb{Q})$ local martingale stopped before $\tau$ 
and $L$ be a $\gg$ predictable process.
The process $L$ is $W$ integrable on $[0,T]$ in the sense of $(\gg,\mathbb{Q})$ local martingales
if and only if $L'$ is $W'$ integrable on $[0,T]$ in the sense of $(\ff,\mathbb{P})$ local martingales.
If so, then (with  ``$\centerdot$ in $(\ff,\mathbb{P})$'', resp.~``$\centerdot$ in $(\gg,\mathbb{Q})$'', in reference to the stochastic integrals in $(\ff,\mathbb{P})$ and $(\gg,\mathbb{Q})$)
$$
\left(L'  \centerdot W' \ \mbox{ in $(\ff,\mathbb{P})$}\right)^{\tau-}
=
\left(L  \centerdot W \ \mbox{ in $(\gg,\mathbb{Q})$} \right) \mbox{ holds on $[0,T]$}.
$$
\ethe

\proof In view of \protect\shortciteN[Definition 9.1]{HeWangYan92}, we only need to check the local integrability of the processes $\sqrt{\int_0^{t} L^2_sd[W,W]_s}$ and $\sqrt{\int_0^t (L')^2_sd[W',W']_s}$ under respectively $(\gg,\mathbb{Q})$ and $(\ff,\mathbb{P})$. But, on $[0,T]$, $[W',W']$ is the $\ff$ optional reduction of $[W,W]$, by Lemma \ref{l:MN}. Hence the local integrabilities above are equivalent, because of Lemma \ref{integrabilityT}. 

To prove the identity between the stochastic integrals when they exist, we first note that the identity holds for any $L$ in the class of $\gg$ predictable bounded step processes. By monotone class theorem, this
is then extended to the class of $\gg$ predictable bounded processes $L$. By stochastic dominated convergence, i.e.~Theorem 9.30 in \shortciteN{HeWangYan92}, this
is extended further to all $\gg$ predictable processes $L$ which are $W$ integrable under $(\gg,\mathbb{Q})$.~\ok

\section{Transfer of Random Measures Stochastic Integrals}
\label{s:rmsi}

Given a Polish space $E$ endowed with 
its Borel $\sigma$ algebra $\cB (E)$\index{b@$\cB (E)$}, we recall from \shortciteN[Theorem 11.13]{HeWangYan92} that, for any 
(optional) integer valued random measure $\pi$, there exists an $E$ valued 
optional process $\beta$ and an 
optional thin set,
of the form
$\cup_{n\in\N}[\theta_n]$
for some sequence of 
stopping times
$(\theta_n)_{n\geq 0}$,
such that
\beql{d:piprel}\pi=\sum_{s} \boldsymbol\delta_{(s,{\beta}_s)}\ind_{\{s\in \cup_{n\in\N}[\theta_n]
\}} \eeql 
(where $\boldsymbol\delta_{(s,{\beta}_s)}$ is the Dirac measure at $(s,{\beta}_s)$).
By definition, for any nonnegative
$\cA \prad\cB(\R_+)\prad \cB (E)$
measurable
function $\Psi,$ 
\beql{d:pi}
\Psi *\pi  = \sum_{s  < \cdot } \Pssi _s(\beta_s)\ind_{\{s\in \cup_{n\in\N}[\theta_n]
\}} = \sum_{\theta_n  < \cdot} \Psi (\theta_n,\beta_{\theta_n})\ind_{\{\theta_n<\infty\}}.
\eeql

\bl The $\gg$ optional integer valued random measure $\pi$
on $\mathbb{R}_+\times  E $ admits
an $\ff$ optional reduction, i.e.
 an $\ff$ optional integer valued random measure $\pi'$ on $\mathbb{R}_+\times  E $ such that $\ind_{[0,\tau)}.\pi=\ind_{[0,\tau)}{\cdot}\pi'$. 
\el
 
\proof  We have, for
any nonnegative
$\cA\prad\cB(\R_+)\prad \cB (E)$
measurable
function $\Psi$,
\bel
& \Psi* (\ind_{[0,\tau)}.\pi)=
\sum_{s\b{< \cdot}} \ind_{\{s<\tau\}} \Pssi _s(\beta'_s)\ind_{\{s\in \cup_{n\in\N}[\theta_n]
\}}\\&\qqq
=
\sum_{s\b{< \cdot}} \ind_{\{s<\tau\}}\Pssi _s({\beta}'_s)\ind_{\{s\in \cup_n[\theta'_n]\}}= \Psi* (\ind_{[0,\tau)}{\cdot}\pi'),
\eel 
where $\pi'=\sum_{s\b{< \cdot}} \boldsymbol\delta_{(s,{\beta}'_s)}\ind_{\{s\in \cup_n[\theta'_n]\}} $
 defines an $\ff$ optional integer valued random measure,
 by \shortciteN[Theorem 11.13]{HeWangYan92}.~\ok \\

\b{In the remainder of the paper}, we fix the space
$E$ and
a $\gg$ optional integer valued random measure  $\pi$, with the related notation
 in the above.
We introduce the spaces of random functions $\widehat{\cP }(\ff)=\cP (\ff)\prad\cB(E)$ and $\widehat{\cP }(\gg)=\cP (\gg)\prad\cB(E)$. 
We denote the $(\ff,\mathbb{P})$ compensator of $\mu=\pi'$ by $\nu.$

\bl\label{proG}
The $(\gg,\mathbb{Q})$ compensator of $\ind_{[0,\tau)}{  \cdot}\mu$ is $\ind_{[0,\tau]}  \cdot\nu$ on $[0,T]$. 
\el

\proof By Lemma \ref{integrabilityT}, 
for any $\Psi\in\widehat{\cP }(\gg)$ such that the process $|\Psi|*\pi$ is $(\gg,\mathbb{Q})$ integrable on $[0,T]$, the processes $|\Psi'|*\mu$ and $|\Psi'|*\nu$ are $(\ff,\mathbb{P})$ locally integrable on $[0,T]$.
It follows that the process
$$
P = \Psi'*\mu - \Psi'*\nu
$$
is an $(\ff,\mathbb{P})$ local martingale on $[0,T]$ (cf.~\shortciteN[p.~301]{HeWangYan92}). By the condition (C), the stopped process 
$$
P^{\tau-} = \ind_{[0,\tau)}\Psi'*\mu - \ind_{[0,\tau)}\Psi'*\nu
= \ind_{[0,\tau)}\Psi*\mu
 - \ind_{[0,\tau]}\Psi*\nu
$$
is a $(\gg,\mathbb{Q})$ local martingale on $[0,T]$, where $\ind_{[0,\tau)}\Psi*\nu=\ind_{[0,\tau]}\Psi*\nu$ because $\tau$ avoids the predictable stopping times. As $\ind_{[0,\tau]}.\nu$ is a $\gg$ predictable random measure, this proves the lemma.~\ok

\bethe\label{itf-randommeasure}
For any $\Psi\in\widehat{\cP }(\gg)$, $\Psi$ is $(\ind_{[0,\tau)}{  \centerdot}\mu - \ind_{[0,\tau]}{  \centerdot}\nu)$ 
stochastically 
integrable
in $(\gg,\mathbb{Q})$ on $[0,T]$ if and only if $\Psi'$ is $(\mu-\nu)$ 
stochastically 
integrable in $(\ff,\mathbb{P})$ on $[0,T]$. If so, then
$$
\left(\Psi'*(\mu-\nu)\ \mbox{ in $(\ff,\mathbb{P})$}\right)^{\tau-}=\left(\Psi*(\ind_{[0,\tau)}.\mu-\ind_{[0,\tau]}.\nu)  \mbox{ in $(\gg,\mathbb{Q})$}  \right)  \mbox{ holds on $[0,T]$}.
$$
\ethe

\proof In view of \citeN[Definition 11.16]{HeWangYan92}, the integrability relationship between $\Psi$ and $\Psi'$ is the consequence of Lemma \ref{integrabilityT}. 
To prove the identity between the corresponding integrals when they exist, we note that $$
(\Psi'*(\mu-\nu))^{\tau-}  \mbox{ and }\
\Psi*(\ind_{[0,\tau)}.\mu-\ind_{[0,\tau]}.\nu)
$$
are $(\gg,\mathbb{Q})$ purely discontinous local martingales. By virtue of ~\citeN[Theorem 7.42 and Definition 11.16]{HeWangYan92}, they are then equal because they have the same jumps, namely
\bel
&\Delta_t(\Psi'*(\mu-\nu))^{\tau-}
= 
\big(\Psi'_t(\beta'_{t})\ind_{\{t\in  \cup_{n\in\N}[\theta'_n]\}}-\int_{\{t\}\times E}\Psi'_s(e)\nu(ds,de)\big)\ind_{\{t<\tau\}}\\&\quad\quad
= 
\big(\Psi_t(\beta_{t})\ind_{\{t\in  \cup_{n\in\N}[\theta_n]\}}-\int_{\{t\}\times E}\Psi_s(e)\nu(ds,de)\big)\ind_{\{t<\tau\}}\\&\quad\quad
=
\Delta_t\big(\Psi*(\ind_{[0,\tau)}.\mu-\ind_{[0,\tau]}.\nu)\big),
\eel
as $\ind_{[0,\tau)}.\nu=\ind_{[0,\tau]}.\nu$ (because $\tau$ avoids the 
predictable stopping times).~\finproof

\section{Transfer of Martingale Representation Properties}\label{s:marepr}

We consider martingale representations with respect to martingales and compensated jump measures as in \citeN{Jacod1979}, which corresponds to the notion of weak representation
in \shortciteN{HeWangYan92}. As in \shortciteN{HeWangYan92}, when no jump measure is involved, we talk of strong representation.

Let $W$ be a $d$ variate $(\gg,\mathbb{Q})$ local martingale stopped before $\tau$ and at $T$.
We assume the random measure $\pi$ of Section \ref{s:rmsi} stopped before $\tau$ and at $T$, in the sense that
$\cup_{n\in\N}[\theta_n]
\subset (0,\tau) \cap (0,T]$.
 We write
$B=(W')^T$, $\mu=(\pi')^T$. Let
$\rho$ and $\nu$
denote the $(\gg,\mathbb{Q})$ compensator of $\pi$ and the $(\ff,\mathbb{P})$ compensator of $\mu$, so that $\rho=\ind_{[0,\tau]}{\cdot}\nu,$ by Lemma \ref{proG}.

\bl\label{l:reprtr} Given $({\cP }(\gg))^{\prad d}$ and $\widehat{\cP }(\gg)$ measurable integrands $L$ and $\Psi$,
if \beql{e:reprG}M=L\centerdot  W  + \Psi * (\pi  -  \rho)\eeql
holds in $(\gg,\Q)$,
then 
$(M')^T=L' \centerdot  B  + \Psi' * (\mu  -  \nu ) $ holds in 
$(\ff,\mathbb{P})$. 
 
Conversely, given $({\cP }(\ff))^{\prad d}$ and $\widehat{\cP }(\ff)$ measurable integrands $K$ and $\Phi$,
if \beql{e:reprF}P=K \centerdot  B  + \Phi * (\mu  -  \nu ) \eeql holds in 
$(\ff,\mathbb{P})$,
then
$P^{\tau-}= K
\centerdot  B^{\tau-}  + \Phi * (\ind_{[0,\tau)}{\cdot}\mu  -  \ind_{[0,\tau]}{\cdot}\nu )$
holds in $(\gg,\Q)$ on $[0,T]$.
\el
\proof This is the consequence of Theorems \ref{itf-usualmartingale} and \ref{itf-randommeasure}. \ok

\brem\label{e:uniq} In the representation \qr{e:reprG}, the integrands $L$ and $\Psi$ corresponding to a given process $M$ are unique modulo $d[W,W]$ 
(with $d[W,W]_s\mbox{-}a.e.$ in the multivariate sense of
\citeN{Jacod2003})
and $\rho$
negligible sets, respectively.
Likewise, in the representation \qr{e:reprF}, the integrands $K$ and $\Phi$ corresponding to a given process $P$ are unique modulo $d[B,B]$ (with $d[B,B]_s\mbox{-}a.e.$ in the multivariate sense of
\citeN{Jacod2003}) and $\nu$ negligible sets.~\finproof
\erem
\noindent
As an immediate consequence of Lemma \ref{l:reprtr}:
\bt\label{mrtcoroll} 
The space $\mathcal{M}_{\tau- \wedge T}(\gg,\mathbb{Q})$ admits a 
weak
representation by $W$ and $\pi$ if and only if the space $\mathcal{M}_T(\ff,\mathbb{P})$ admits a 
weak
representation by $B=W'$ and $\mu=\pi'$.~\finproof
\et

\noindent
Applying Theorem \ref{mrtcoroll} with $\mu\equiv 0$, one obtains the strong martingale representation transfer property.

We refer the reader to \citeN{gapeev2021projections,gapeev2022projections} for other transfers of martingale representation properties, in respective Brownian and marked point process enlargement of filtration setups (progressive but also initial as already before in \citeN{fontana2018strong}) satisfying Jacod's equivalence hypothesis, i.e. the existence of a  positive  $\ff$ conditional density for $\tau$,
as opposed to a semimartingale progressive enlargement of filtration setup under the condition (C) in this work. See also \citeN{jeanblanc2015martingale} or (until $\tau$) \citeN{choulli2020martingale} and  (also after $\tau$)  \citeN{choulli2022representation} for 
rather general transfers of martingale representation properties in a progressive enlargement of filtration setup. From a technical viewpoint, our setup stopped before $\tau$ 
is elementary
once the underlying Theorems \ref{itf-usualmartingale} and \ref{itf-randommeasure} are in place. 
Of
course one cannot say anything beyond $\tau$ in our setup, but our motivating application of Example \ref{e:cy} never requires to go beyond $\tau$.

\section{Semimartingale Characteristic Triplets Transfer Formula}
{
\def\d{\boldsymbol{\alpha}}\def\d{b}
\def\a{\boldsymbol{\sigma}}\def\a{a}
\def\c{\boldsymbol{\nu}}\def\c{c}
\def\be{\boldsymbol{\gamma}}\def\be{\beta}
\def\pip{\pi'}\def\pip{\pi^{(X')^T}} 

Let there be given a semimartingale $X$ stopped before $\tau$
in some filtration $\hh $ under a probability measure $\mathbb{M}$, with jump measure ${\pi^X}$. The characteristic triplet 
of $X$ is composed of:
\bel&
\d ^{X,\hh , \mathbb{M}},  \mbox{the drift
part of the truncated
semimartingale $X - (x\ind_{\{|x|>1\}})_{*}{\pi^X}$};\\ 
&\a ^{X,\hh , \mathbb{M}}= \mbox{ $\croc{X^c,X^c}$, the angle bracket of the continuous martingale part of $X$}\\
& \mbox{(i.e. the diffusion part of $X$);}\\ 
&\c ^{X,\hh , \mathbb{M}}= \mbox{ $({\pi^X})^{p, \hh , \mathbb{M}}$, the predictable dual projection of ${\pi^X}$,}\\
& \mbox{called in \shortciteN{HeWangYan92} the L\'evy system of $X$ (i.e.~the} \\
& \mbox{extension to a semimartingale setup of the notion of a L\'evy measure)}.
\eel
 
The following results show (essentially) that the $(\ff, \mathbb{P})$ characteristic triplet of the optional reduction $X'$ of 
a $(\gg, \mathbb{Q})$ semimartingale stopped before $\tau$, $X$, is the predictable reduction of the $(\gg, \mathbb{Q})$ characteristic triplet of $X$. Moreover, if the $(\gg, \mathbb{Q})$ semimartingale $X$ is special, then so is $X'$ and
the $(\ff, \mathbb{P})$ drift of $X'$ is the predictable reduction of the $(\gg, \mathbb{Q})$ drift of $X$.

\bethe\label{tripletreduction}
Let $X=X^{\tau-}$ be a $(\gg, \mathbb{Q})$ semimartingale stopped before $\tau$ and at $T$.
We have
\beql{e:tri}
\left(
\d ^{X,\gg, \mathbb{Q}},
\a ^{X,\gg, \mathbb{Q}},
\c ^{X,\gg, \mathbb{Q}}
\right)
=
\left(
(\d ^{(X')^T ,\ff, \mathbb{P}})^{\tau},
(\a ^{(X')^T ,\ff, \mathbb{P}})^{\tau},
\ind_{[0,\tau]}.\c ^{(X')^T ,\ff, \mathbb{P}}
\right) .
\eeql 
\ethe

\proof
First note that $(X')^T$ is an $(\ff, \mathbb{Q})$ semimartingale,  
by \citet[Lemma 
 6.5]{Song16a} (with $\mathtt{S}_T>0$ a.s.~under the condition (C)), hence an  $(\ff, \mathbb{P})$ semimartingale as well,  which justifies the writings in the right hand side of \eqref{e:tri}.
Denoting by $\pip $ the $(\ff, \mathbb{P})$ jump measure of $(X')^T$,
we have  
$
\ind_{[0,\tau)}.{\pi^X}
=
\ind_{[0,\tau)}.\pip .
$ 
So
\begin{equation}\label{truncmart}
X - (x\ind_{\{|x|>1\}})_{*}{\pi^X}
=
\left(
(X')^T  - (x\ind_{\{|x|>1\}})_{*}\pip 
\right)^{\tau-}
=
(P^c)^{\tau-}+(P^d)^{\tau-}+(\d ^{(X')^T ,\ff, \mathbb{P}})^{\tau-},
\end{equation}
where $P$ is the $(\ff, \mathbb{P})$ canonical Doob--Meyer local martingale component of the $(\ff,\P)$ special semimartingale $(X')^T  - (x\ind_{\{|x|>1\}})_{*}\b{\pip }$, with continuous and purely discontinuous parts $P^c$ and $P^d$. By the condition (C), $P^{\tau-}$ is a $(\gg, \mathbb{Q})$ local martingale. Therefore, we conclude from (\ref{truncmart}) that$$
\d ^{X,\gg, \mathbb{Q}}
=
(\d ^{(X')^T ,\ff, \mathbb{P}})^{\tau-}
=
(\d ^{(X')^T ,\ff, \mathbb{P}})^{\tau} 
$$
(as $\Delta_{\tau}\d ^{(X')^T ,\ff, \mathbb{P}}=0$, because $\tau$ is totally inaccessible.) Now, applying Lemma \ref{proG} with $E=\R$, we also conclude$$
\c ^{X,\gg, \mathbb{Q}}
=
({\pi^X})^{p,\gg, \mathbb{Q}}
=
(\ind_{[0,\tau)}.{\pi^X})^{p,\gg, \mathbb{Q}}
=
\ind_{[0,\tau]}.(\pip)^{p,\ff, \mathbb{P}}
=
\ind_{[0,\tau]}.\c ^{(X')^T ,\ff, \mathbb{P}}.
$$
Finally, according to the second and third bijections in \qr{e:bijdtc}, we have
 $$(P^c)^{\tau-}\in\mathcal{M}^{c}_{\tau  \wedge T}(\gg,\mathbb{Q})\sp (P^d)^{\tau-}\in\mathcal{M}^{d}_{\tau- \wedge T}(\gg,\mathbb{Q}).$$
 Hence
we conclude from (\ref{truncmart}) that $
X^c = (P^c)^{\tau-}
$
is the continuous local martingale part of $X$ in $(\gg, \mathbb{Q})$ and therefore $$
\a ^{X,\gg, \mathbb{Q}}
=
[X^c,X^c]
=
[(P^c)^{\tau-}, (P^c)^{\tau-}]
 =
[ P^c ,  P^c ]^{\tau} 
=
(\a ^{(X')^T ,\ff, \mathbb{P}})^{\tau}.\ \ok
$$
\bcor
Suppose that a $(\gg,\Q)$ semimartingale $X=X^{\tau-}$ is special on $[0,T]$. Then $(X')^T $ is an $(\ff,\P)$ special semimartingale. Denoting by $\be ^{X,\gg, \mathbb{Q}}$ and $\be ^{(X')^T ,\ff, \mathbb{P}}$ the $(\gg, \mathbb{Q})$ drift of $X$ and the $(\ff, \mathbb{P})$ drift of $(X')^T ,$ we have
\beql{e:dri}
\be ^{X,\gg, \mathbb{Q}}
=
(\be ^{(X')^T ,\ff, \mathbb{P}})^\tau.
\eeql 
\ecor

\proof As \b{$(X')^T $ is already known to be an $(\ff,\P)$ semimartingale}
and because a special semimartingale means one with locally integrable jumps,
the special feature of $(X')^T $ follows from Lemma \ref{integrabilityT}.
Note that, by \citep*[Lemma 7.16 and Theorem 11.24]{HeWangYan92}, the function $|x|\ind_{\{|x|>1\}}$ is $\c ^{(X')^T ,\ff, \mathbb{P}}$ integrable on $[0,T]$. Consequently$$
\be ^{(X')^T ,\ff, \mathbb{P}}
=
\d ^{(X')^T ,\ff, \mathbb{P}}
+
(x\ind_{\{|x|>1\}})_{*}\c ^{(X')^T ,\ff, \mathbb{P}}.
$$
The analogous $(\gg,\Q)$ relationship holds for $X$. Hence \qr{e:dri} follows from \qr{e:tri}. \ok
 }

\section{Markov 
Transfer Formulas}\label{s:markov} 

In this section we study the transfer of 
Markov properties between $(\gg,\Q)$ and $(\ff,\P)$.
The reader is referred to  \citeN[proof of Proposition (60.2)]{Sharpe88}
regarding the definition of the semigroup generated by a Markov family.

We suppose 
that the filtration
$\gg$ is generated by a  \b{$(\gg,\Q)$}
quasi-left continuous strong Markov semimartingale $\Xb $  with state space
 $\mathbb{R}^d$.
We denote by an index ${\cdot}^{(t)}$ \b{everything related to the Markov process $X$ translated by time $t$}.
We assume that $\tau$ is a terminal time of $X$, i.e.~(see e.g.~\citeN[(3.7) Remark p.108]{BlumenthalGetoor2007}\footnote{or
Definition 1.3 page 98 in the original 1968 edition of their book.})
$
\tau=\b{\tau^{(t)}}
+t \mbox{ if }\tau>t $\footnote{or  ``$\tau=\tau\circ \theta_t +t$'' in the classical notation where $\theta$ is the translation operator associated with the Markov process $X$ (but we already use $\theta$ to denote a stopping time in this paper).}.
We assume further that the
$(\gg,\Q)$ intensity process of $\tau$ takes the form $\gamma(\Xb _{\cdot})\ind_{(0,\tau]},$ for some $\cB(\R^d )$ measurable function $\gamma \ge 0$. 
Let $$
\mathtt{M}_s=e^{\Gamma_s}\ind_{\{s<\tau\}}
=e^{\int_0^{s}\gamma(X_u)du}\ind_{\{s<\tau\}}.
$$

\bl
$\mathtt{M}$ is a multiplicative functional of $X$, i.e. $\mathtt{M}_t=\mathtt{M}_s  \mathtt{M}_{t-s}^{(s)}  $ holds for any $t\ge s\ge 0$, 
and a $(\gg,\mathbb{Q})$ local martingale. The multiplicative functional $M$ defines a probability transition function $(\mathcal{T}_t)_{t\in\mathbb{R}_+}$.
\el

\proof The first part can be checked by definition of $\mathtt{M}$. For the second part, we check by the Dol\'eans-Dade exponential formula that $$
\mathtt{M} =\mathcal{E}(-\ind_{[\tau,\infty)}+\Gamma^{\tau}). 
$$ 
The last part follows from \b{\citeN[(65.3), proof of Proposition(56.5)]{Sharpe88}, \citeN[Hypothesis (62.9)]{Sharpe88} and \citeN[Theorem (62.19)]{Sharpe88}}. \ok

\bethe\label{th:strm}
The reduction $X'$ of $X$ is an $(\ff,\mathbb{P})$ strong Markov process on $[0,T], $
with the transition semigroup $(\mathcal{T}_t)_{t\in[0,T]}$.
\ethe

\proof
For $A\in\mathfrak{F}_s,$ $h$ Borel bounded and $0<s<s+t\leq T$, we have 
\bel
&\mathbb{E}[\ind_A h (X_{t+s})e^{\int_0^{t+s}\gamma(X_u)du}\ind_{\{t+s<\tau\}}]\\
&\qq=\mathbb{E}[\ind_Ae^{\int_0^{s}\gamma(X_u)du}\ind_{\{s<\tau\}}
\b{ h (X_{t}^{(s)}) e^{\int_0^{t}
\gamma(X_u^{(s)})du}\ind_{\{t<\tau^{(s)}\}} }]\\
&\qq=
\mathbb{E}[\ind_Ae^{\int_0^{s}\gamma(X_u)du}\ind_{\{s<\tau\}}
\mathcal{T}_t h (X_{s})],
\eel
which is rewritten in terms of $X'$ through the first expectation transfer formula
\qr{t:etf}
as
$$
\dcb
&&\Ep[\ind_A h (X'_{t+s})]
=
\Ep[\ind_A (\mathcal{T}_t h) (X'_{s})].
\dce
$$
This proves that $X'$ is an $(\ff,\mathbb{P})$ Markov process with the transition semigroup $(\mathcal{T}_t)_{t\in[0,T]}$. If we rewrite the above computation for $s$ replaced by an $\ff$ stopping time $\sigma$, we prove that $X'$ is an $(\ff,\mathbb{P})$ strong Markov process on $[0,T]$. \ok\\

The next question is how to determine the generator of the semigroup $(\mathcal{T}_t)_{t\in[0,T]}$.
We suppose that 
the Markov process $X$ is of the form $X=(Y,Z),$ \b{for some 
process $Y$ stopped before $\tau$ and some 
process $Z$ constant ($=0$, say) before $\tau$, the role of which is to store some information, encoded into the jump of $Z$ at $\tau$, about what happens at $\tau$;
\bex\label{e:dc} We may consider for $X$ the following dynamic copula models of portfolio credit risk, with any of the portfolio default times in the role of $\tau$ in this paper:
{\rm\hfill\break \textbf{(i)}}  the dynamic Marshall-Olkin copula (DMO) model, 
shown in \citeN[Theorem 9.2]{CrepeySong15FS} to satisfy the condition (C) for $\P=\Q$ there {(case of a strict pseudo-stopping time in the terminology of \citeN[Definition 2.1]{jeanblanc2020characteristics})};
{\rm\hfill\break \textbf{(ii)}} the  
dynamic Gaussian copula (DGC) model,
shown in \citeN{CrepeySong16b} to satisfy the condition (C) with $\P\neq\Q$,  provided the
correlation coefficient $\varrho>0$ in the model is small enough. In particular, the condition (C) holds in the univariate DGC model (there is then no correlation $\varrho$ involved) of \sr{ss:ex} below.
\eex}
\noindent
Suppose further that $X$ solves the following  $(\gg,\mathbb{Q})$ martingale problem:$$
 v(X_t) - \int_0^t \cL v(X_s) ds \mbox{ is a $(\gg,\mathbb{Q})$ local martingale for all $v\in\mathcal{D}(\cL)$,}
$$ 
where $\cL$ is \b{the generator of $X$, with domain $\mathcal{D}(\cL)\subseteq$ the set of the $\cB(\R^d)$ measurable bounded functions.}  
For  $u\equiv u(y)$ we define $\hat{u}\equiv\hat{u}(y,z)$ by  $\hat{u}(y,z)=u(y)$.
Let 
$
\mathcal{D}'
=
\{u \equiv u(y); \hat{u}\in \mathcal{D}(\cL)\}$ and let $\cL'$ be the operator
on $\mathcal{D}'$ defined by
\beql{e:gens}
&(\cL' u)(y)
=
(\cL \hat{u} ) (y,0)\sp  u\in \mathcal{D}'  .\eeql
\bethe\label{th:gene}
We suppose that $(\mathcal{D}',\cL')$ satisfies the conditions of \citeN[Theorem 4.1 of Chapter 4, p.182]{EKu}. 
Then $X'=(Y',0)$, $Y'$ is an $(\ff,\mathbb{P})$ strong Markov process on $[0,T]$, and the generator of $Y'$ is 
\b{an
 extension of} $(\mathcal{D}',\cL')$.
\ethe
\proof
Clearly, $X'=(Y',0)$. Hence, Theorem \ref{th:strm} implies that $Y'$ is an $(\ff,\mathbb{P})$ strong Markov process on $[0,T]$. For $u\in\mathcal{D}'$,$$
\dcb
&& \hat{u}(X_t) - \int_0^t \cL \hat{u}(X_s) ds =
 u(Y_t) - \int_0^t \cL \hat{u}(Y_s, 0) ds \dce$$
is a $(\gg,\mathbb{Q})$ local martingale. 
As $Y$ is stopped before $\tau$, 
$$
 u(Y_t)^{\tau-} - \int_0^{t\wedge \tau} \cL \hat{u}(Y_s, 0) ds  \mbox{ is a $(\gg,\mathbb{Q})$ local martingale stopped before $\tau.$} 
$$
The first bijection in \eqref{e:bijdtc} then implies that
$$
 u((Y')^T_t) - \int_0^{t} \cL' u((Y')^T_s) ds  \mbox{ is an $(\ff,\mathbb{P})$ local martingale on $[0,T]$}.
$$ 
Therefore, $Y'$ is the solution of the $(\ff,\mathbb{P})$ martingale problem associated with $\cL'$ on $[0,T]$. The result then follows from \citeN[Theorem 4.1 of Chapter 4, p.182]{EKu}.~\ok

\section{BSDE Transfer Properties} \label{s:bsde}
 
In this section,
$\tau$ satisfying the condition (C) on $[0,T]$ as before, we reduce a $(\gg,\Q)$ backward stochastic differential equation (BSDE) stopped before 
$\tau$ and at $T$ to a simpler $(\ff,\P)$ BSDE stopped at $T$.

We suppose $E$ Euclidean and $(E,\cB(E))$ endowed
with a $\sigma$ finite measure $m$ integrating ($1 \wedge |e|^2$) on $E$.
We consider
 the space $\mathbb{L}_0 $ of the $\cB(E)$ measurable functions  $u$ endowed with the topology of convergence in measure induced by $m$.

Given  a $\cP (\gg)\prad \cB(\mathbb{R})\prad \cB(\mathbb{R}^d) \prad {{\cB} (\mathbb{L}_0   )
}	$
 measurable function $g=g_t(z,l,\psi)$,
we can define, {by monotone class theorem},
a $\cP (\ff)\prad \cB(\mathbb{R})\prad \cB(\mathbb{R}^d) \prad {{\cB} (\mathbb{L}_0   )
}	$
reduction
$g'=g'_{t}(z,l,\psi)$ of $g$ such that $\ind_{(0,\tau]}g =\ind_{(0,\tau]} g'$.
Let $A$ be a $\gg$
finite variation (\cadlag) 
process.

Adopting the setup of \sr{s:marepr}, we consider the $(\gg,\Q)$ BSDE with data $(g,A)$ and solution sought for as a triplet  
$(Z,L,\Psi)$, where $Z$ is
 a
 $\gg$ adapted process, 
  $L$ is
a $(\cP(\gg))^{\prad d}$  measurable process
integrable against
$B^{\tau-}$ in $(\gg,\Q)$,
and $\Psi$ is a $\hat{\cP}(\gg)$ measurable function 
 stochastically integrable against $(\ind_{[0,\tau)}{\cdot}\mu  -  \ind_{[0,\tau]}{\cdot}\nu )$
in $(\gg,\Q)$,
 satisfying in $(\gg,\Q)$
\beql{e:Z}
\left\{
\dcb
\int_0^{\tau\wedge T}|g_s(Z_{s-}, L_s,  \Psi_{s})| ds <\infty \mbox{ and }  \\
\qqq \int_0^{\cdot}\indi{s<\tau}|dA_s|\mbox{  is  } 
(\gg,\Q)
\mbox{ locally integrable on } [0,T],\\
\\
Z^{\tau-\wedge T}_t+\int_0^{t\wedge \tau\wedge T}\big(g_s(Z_{s-}, L_s, \Psi_{s}) ds
+
dA^{\tau-}_s\big) \\
\qqq=L\centerdot  B^{\tau-}_t + \Psi * (\ind_{[0,\tau)}{\cdot}\mu  -  \ind_{[0,\tau]}{\cdot}\nu )_t\sp  t\in\R_+ ,\\
\\
\b{Z \mbox{ vanishes on } [\tau\wedge T,+\infty).} 
\dce
\right.\eeql  
We also consider the
$(\ff,\P)$ BSDE with data $(g',A')$ and solution 
 sought for as a triplet  
$(U,K,\Phi)$, where $U$ is
 an $\ff$ adapted process, $K$ is a $(\cP(\ff))^{\prad d}$  measurable process 
integrable against
$B$ in $(\ff,\P)$ on $[0,T]$,
and $\Phi$ is a $\hat{\cP}(\ff)$ measurable function 
 stochastically integrable against $(\mu  - \nu )$
in $(\ff,\P)$ on $[0,T]$,
  satisfying in $(\ff,\P)$
\begin{equation}\label{Fmartpb}
\left\{
\dcb
\int_0^T | g'_s (U_{s-}, K_s, \Phi_s)|ds<\infty  
\mbox{ and }  \\
\qqq\int_0^{\cdot}|dA'_s| \mbox{  is  } 
(\ff,\P)
\mbox{ locally integrable on } [0,T],
\\\\
U^T_t+\int_0^{t\wedge T}\big(g'_s(U_{s-}, K_s, \Phi_s) ds+dA'_s \big)
=K\centerdot  B_t + \Phi * (\mu  -  \nu )_t\sp  t\in\R_+ ,
\\
\\
\b{U^T \mbox{ vanishes on } [ T,+\infty) \mbox{ (i.e. $U_T=0$)}.} 
\dce
\right.
\end{equation}
Note that the $(\gg,\Q)$ BSDE \qr{e:Z} is stopped at $\tau-\wedge T$ \b{(in particular, the terminal condition $Z_T=0$ only holds on $\{T<\tau\}$)},
whereas the $(\ff,\P)$ BSDE \qr{Fmartpb} is stopped at $T.$

\bex \label{e:cy}
Given a bank with  default time $\tau$, a $\gg$ stopping time $\theta$ representing the default time of a client of the bank, and a nonnegative $\gg$ optional process $G$ representing the liability of the client to the bank, then
the process $A=\int_0^{\cdot} G_s \boldsymbol\delta_{\theta}(ds)$ represents the counterparty credit exposure of the bank to its client. In this case
$$|dA_s|= G_s \boldsymbol\delta_{\theta}(ds)\sp
A^{\tau-}=\int_0^{\cdot} \indi{s<\tau} G_s \boldsymbol\delta_{\theta}(ds)\sp
A'=\int_0^{\cdot}{ G'_s }\boldsymbol\delta_{\theta'}(ds) .$$
The coefficient $g$ represents the risky funding costs  of the bank entailed by its own credit riskiness.
For the reason explained in the second paragraph of \sr{s:intro}, 
all cash flows are stopped before the bank default time $\tau$. 
This results in a BSDE of the form \qr{e:Z} for the valuation of counterparty risk (CVA) and of its funding implications to the bank (FVA). The cost of capital (KVA) also obeys an equation of the form \qr{e:Z}: see \citeN[Eqns. (2.12), (2.13), and (2.17)]{Crepey21}.
\eex

\subsection{Transfer of Local Martingale  Solutions}

The result that follows states the equivalence between 
the $(\gg,\Q)$ BSDE \qr{e:Z} and the $(\ff,\P)$ BSDE \qr{Fmartpb}
considered within the above-introduced spaces of solutions 
for the triplets $(Z,L,\Psi)$ and $(U,K,\Phi)$, called local martingale solutions henceforth (in reference to the fact that the right-hand sides in the second lines of \qr{e:Z} and \qr{Fmartpb} are then respectively $(\gg,\Q)$ and $(\ff,\P)$ local martingales).
 
\bt\label{t:bsdeequiv} The $(\gg,\Q)$ BSDE \qr{e:Z} and the $(\ff,\P)$ BSDE \qr{Fmartpb} are equivalent in their respective spaces of local martingale solutions.
Specifically, if $(Z,L,\Psi)$ solves (\ref{e:Z})
in
$(\gg,\mathbb{Q})$, then 
$(U,K,\Phi)=(Z,L,\Psi)'$ solves \qr{Fmartpb} in $(\ff,\mathbb{P})$.  
Conversely, if $(U,K,\Phi)$ solves \qr{Fmartpb} in $(\ff,\mathbb{P})$, then
  $(Z,L,\Psi)=(\b{ \ind_{[0,\tau)}U},
\ind_{[0,\tau]}K,\ind_{[0,\tau]}\Phi)$ solves \qr{e:Z} in $(\gg,\mathbb{Q})$.  
\et

\proof Through the correspondence stated in the theorem between the involved processes:
\begin{itemize}
\item
The equivalence between the Lebesgue integrability conditions (first lines) in \qr{e:Z}
and \qr{Fmartpb}
follows from Lemma \ref{integrabilityT}; 
\item
The equivalence between the martingale conditions (second lines) in \qr{e:Z} 
and \qr{Fmartpb} follows from Theorems \ref{itf-usualmartingale} and \ref{itf-randommeasure}; 
\item
The terminal condition  for $U$
in \qr{Fmartpb} obviously implies the one \b{for $Z=\ind_{[0,\tau)}
U$} in  \qr{e:Z}, whereas the terminal condition in (\ref{e:Z}) implies
$Z_{T}\indi{T<\tau}=0$, hence
by taking the $\mathfrak{F}_{T}$ conditional expectation:
$$0=\mathbb{E}[Z_{T}\indi{T<\tau}|\mathfrak{F}_{T }]=\mathbb{E}[Z'_{T}\indi{T<\tau}|\mathfrak{F}_{T }]=Z'_{T} \ttS_T,$$
yielding $U_T=Z'_{T}=0$ (as $\ttS_T$ is positive under the condition (C)).~\finproof 
\end{itemize}

\subsection{Transfer of 
Square
Integrable Solutions}

We now consider
the $(\gg,\Q)$ BSDE \qr{e:Z} and the $(\ff,\P)$ BSDE \qr{Fmartpb}
within suitable spaces of square integrable solutions.
 
We assume that the compensator $\nu$ of $\mu=\pi'$ is given as
$\zeta_t (e) m (de)dt,$
where $\zeta$ is a nonnegative and bounded integrand in $\cP(\ff)$. 
We write, for any $t\ge 0$ 
and $\cB(E)$ measurable function $u$,
$$  
 |u |_t^2= \int_E u(e)^2 \zeta_t (e) m (de).$$
Let also $Y^*_t =\sup_{s\in[0,t]} |Y_s|$, for any \cadlag process $Y$.
\bl \label{l:norms}
For any
real valued c\`adl\`ag $\ff$ 
adapted process 
$V$,  respectively nonnegative $\ff$ predictable process $X$,
we have
\begin{eqnarray}\label{e:norms1}&&
 \mathbb{E}\Big| V_0^2+ \int_0^T \ e^{\int_0^s\gamma_u du} \ind_{\{s<\tau\}}d(V^*)_s^2\Big]    = 
\Et[(V^*)_T^2] ;\\
\label{e:norms2}&&\mathbb{E}\Big[\int_0^T \ e^{\int_0^s\gamma_u du} \ind_{\{s<\tau\}}\    
 X_s 
ds \Big] = \Et\Big[ \int_0^T \ 
 X_s  
\ ds \Big] . 
\end{eqnarray}\el
\proof
The formula 
\qr{t:etfpp} used at $t=0$ yields: 
\begin{itemize}  
\item For $A=\int_0^{\cdot} \ e^{\int_0^s\gamma_u du}  d(V^*)_s^2,$
\bel&
 \mathbb{E}\Big[\int_0^T \ e^{\int_0^s\gamma_u du} \ind_{\{s<\tau\}}d(V^*)_s^2\Big]   = 
\Et[(V^*)_T^2] -  \mathbb{E}'[V_0^2] ;
\eel
\item For $A=\int_0^{\cdot} \ e^{\int_0^s\gamma_u du}   X_s 
\ ds $,
\bel 
&\mathbb{E}\Big[\int_0^T \ e^{\int_0^s\gamma_u du} \ind_{\{s<\tau\}}\    
X_s  
ds \Big] = \Et\Big[ \int_0^T \ 
X_s 
\ ds \Big] . ~\finproof
\eel
\end{itemize}

Considering the $(\gg,\Q)$ BSDE \qr{e:Z} for $(Z,L,\Psi)$
and the reduced $(\ff,\P)$ BSDE \qr{Fmartpb} for $(U,K,\Phi)$,
with local martingale solutions (if any) such that 
\beql{e:lms}&(U,K,\Phi)=(Z,L,\Psi)'\sp (Z,L,\Psi)=( \b{\ind_{[0,\tau)}
U},
\ind_{[0,\tau]}K,\ind_{[0,\tau]}\Phi)\eeql
(cf.~Theorem \ref{t:bsdeequiv}), we define 
\bel&
\| (Z,L,\Psi)\|_{2}^2
=\mathbb{E}\Big[|Z_0|^2+\int_0^T \ e^{\int_0^s\gamma_u du} \ind_{\{s<\tau\}}
d(Z^*)_s^2
\Big]\\&\qqq
+
\mathbb{E}\Big[\int_0^T \ e^{\int_0^s\gamma_u du} \ind_{\{s<\tau\}}\ 
\big( |L_s |^2+ |\Psi_s |_s^2\big)
ds
\Big],\\
& 
( \|(U,K,\Phi)\|'_{2} )^2
=
\Et[(U^*)_T^2 ]
+
\Et\Big[
\int_0^T \ 
  \big( |K_s |^2+ |\Phi_s |_s^2\big)
\ ds
\Big] .
\eel

We consider 
the respective subspaces of square integrable solutions of the $(\gg,\Q)$ BSDE \qr{e:Z} and of the $(\ff,\P)$ BSDE \qr{Fmartpb} defined by
$\|\cdot \|_{2}<+\infty$ and \b{$Z=0$ on $[\tau\wedge T, +\infty)$}, respectively
$\|\cdot\|'_{2}<+\infty$ and $U^T=0$ on $[T, +\infty)$ (i.e. $U_T=0$), dubbed $\|\cdot \|_{2}$ and $\|\cdot \|'_{2}$ solutions hereafter.

\bt\label{c:eqp}\label{l:eqnorm} Given local martingale solutions $(Z,L,\Psi)$ to
the $(\gg,\Q)$ BSDE \qr{e:Z}  
and $(U,K,\Phi)$ to the reduced $(\ff,\P)$ BSDE \qr{Fmartpb},
we have 
\beql{e:nt}
\| (Z,L,\Psi) \|_{2}=\| (U,K,\Phi) \|'_{2}
.\eeql 
The $(\gg,\Q)$ BSDE \qr{e:Z} considered in terms of $\|\cdot \|_{2}$ solutions and the $(\ff,\P)$ BSDE \qr{Fmartpb} considered in terms of $\|\cdot \|'_{2}$ solutions are equivalent through the correspondence \qr{e:lms}.
\et

\proof Given respective local martingale solutions $(Z,L,\Psi)$ and $(U,K,\Phi)$ to \qr{e:Z} and \qr{Fmartpb}, then
related through
\qr{e:lms} as seen in Theorem \ref{t:bsdeequiv},
Lemma \ref{l:norms}  applied to $V=U$ and  $X = |K_{\cdot} |^2+ |\Phi_{\cdot} |_{\cdot}^2$ 
proves the transfer of norms formula \qr{e:nt}. 
Given the equivalence of
Theorem \ref{t:bsdeequiv} between \qr{e:Z} and \qr{Fmartpb} in the sense of local martingale solutions, their equivalence in the sense of square integrable solutions follows from \qr{e:nt}.~\finproof

\subsection{Application}

Assuming $\int_0^{T}{|dA'_s|}$  integrable
under $\P$ and a (weak) martingale representation of the form studied in Theorem \ref{mrtcoroll},
we define the processes $R$ and $P$ given as
\bel
&R_t=\Et\big[ \int_{t\wedge T}^Tt dA'_s
\,|\,\cF_t \big]\mbox{ and } P_t= \Et \big[ \int_0^T dA'_s
\,|\,\cF_t \big]\sp t\in\R_+.
\eel
Let $f_s(v,k,\phi)=g'_s( R_{s-}+v ,K^P_s+k ,\Phi^P_s+\phi)$,
where $K^P$ and $\Phi^P$ are the integrands in the representation \qr{e:reprF} of the $(\ff,\P)$ martingale $P$ (cf.~Remark
\ref{e:uniq}).

\bp\label{t:bsdetranf}
Suppose that  $\int_0^{T}{|dA'_s|}$ is $\P$  square  integrable
and
 \begin{itemize}
 \item[{(i)}] the functions $\R\ni v\mapsto f_t(v,k,\phi)\in\R$ are continuous, for each $(k,\phi)\in\R^d\times \mathbb{L}_0$. Moreover, $f$ is monotonous with respect to $v$, i.e.
 $$ (f_t(v_1,k,\phi )- f_t(v_2, k,\phi ))( v_1-v_2)\leq C (v_1-v_2)^2 ;$$ 
  \item[{(ii)}] $\Et\int_0^T \underset{|v|\leq c}{\sup} |f_t(v,0,0 )-f_t(0,0,0 )|dt <\infty$ holds for every positive $c$; 
  \item[{(iii)}] $f$ is Lipschitz continuous with respect to $k$ and $\phi$, i.e.
$$|f_t(v,k_1,\phi_1)-f_t(v,k_2,\phi_2)| \leq C (|k_1-k_2|+ |\phi_1-\phi_2 |_t  ); $$
    \item[{(iv)}] $\Et \int_0^T |f_t(0,0,0)|^2 dt < +\infty.$
 \end{itemize}
Then the $(\gg,\Q)$ BSDE \qr{e:Z} and the $(\ff,\P)$ BSDE \qr{Fmartpb} have unique
 $\|\cdot \|_{2}$ and $\|\cdot \|'_{2}$ solutions, respectively, and these solutions are related through \qr{e:lms}. 
\ep

\proof \b{Note that  $\int_0^{T}{|dA'_s|}$ being  $\P$  square  integrable 
implies that $\E'\big[(R^*)^2_T\big] <\infty$.}
Through the correspondence
\bel
&U=R +V\sp K^U=K^P +K^V\sp \Phi^U=\Phi^P +\Phi^V,\eel 
the 
$(\ff,\P)$ BSDE \qr{Fmartpb} for 
$(U,K^U,\Phi^U)$ 
is equivalent (in both senses of $(\ff,\P)$ local martingale solutions
and $\|\cdot \|'_{2}$ solutions)
to the following 
$(\ff,\P)$ BSDE for 
$(V,K^V,\Phi^V)$:
\begin{equation}\label{e:V}
\left\{
\dcb
 \int_0^T |{f}_s( V_{s-},  K^V_s, \Phi^V_s)|ds   <\infty ,\\\\
V_t^T +\int_0^{t\wedge T} {f}_s( V_{s-},  K^V_s, \Phi^V_s) ds 
=K^V\centerdot  B_t + \Phi^V * (\mu  -  \nu )_t,
\\
\\
V_T=0.
\dce
\right.
\end{equation}
Under the assumptions of the proposition, 
the $(\ff,\P)$ BSDE \qr{e:V} for $(V,K^V,\Phi^V)$ satisfies the
assumptions of \citeN[Theorem 1]{KrusePopier14}.
Hence it
has a unique $\|\cdot \|'_{2}$ solution.
So has
in turn the  
BSDE \qr{Fmartpb}. The result then follows by an application of Theorem \ref{c:eqp}.
\finproof\\

\noindent
\brem \citeN[Theorem 1]{KrusePopier14} is only derived in the case a Poisson measure $\pi$, 
but one can readily check that all their computations and results derived under square integrable assumptions are still valid in our more general integer valued random measure setup. Also, in view of \citeN{bouchard2018unified}, \citeN{KrusePopier14}'s condition of a quasi-left continuous filtration is in fact not needed.
\erem

\noindent
The reader is referred to
\citeN[Lemma B.1, Proposition B.1 and Theorem 6.1]{Crepey21} and \citeN[Section 6]{CrepeyElie16}
for variations on the above results, in the respective cases where $f$ only depends on $v$ (in the notation of Proposition \ref{t:bsdetranf} above) and no martingale representation property needs to be assumed,
or where $f$  is assumed to be Lipschitz (versus monotonous in Proposition \ref{t:bsdetranf}) but also exhibits a dependence on a conditional expected shortfall of a future increment of the martingale part of the solution. 

\brem 
Earlier occurrences of such results are \citeN{Crepey13a1,CrepeySong15FS}, with the difference that these earlier works were about BSDEs stopped at a random time. The more recent papers, instead, with the motivation recalled in the second paragraph of Section \ref{s:intro}, are about BSDEs stopped before a random time: compare e.g. stopping at $\vartheta$ in the second line of \citep[Eqn. (2.1)]{Crepey13a1} versus stopping before $\tau$ in the second part of \eqref{e:Z}.
\erem

{Analogous techniques could be used to simplify $(\gg,\Q)$ optimal stopping or stochastic control problems into reduced  $(\ff,\P)$ reformulations: cf., in the case of BSDEs or control problems stopped at time $\tau$ (as opposed to stopped before $\tau$ in our setup, and without the flexibility induced by our measure change from $\Q$ to a possibly different $\P$), \citeN{KharroubiLim11}
and \citeN{jiao2013optimal} (assuming that a driving $(\F,\Q)$ Brownian motion, stopped at $\tau$, is a $(\G,\Q)$ martingale), \citeN{aksamit2021generalized} (who provide some comparative comments with our approach in their Remark 8.2), or \citeN{alsheyab2021reflected}.}

\section{Conclusion}\label{s:concl}

The present paper provides new developments on the concept of invariance times $\tau$ introduced in \citeN{CrepeySong15c}. More precisely, this paper provides more complete results under stronger (practical) assumptions, summarized as the condition (C), whereas \citeN{CrepeySong15c} was mainly about the conditions (B) and (A) corresponding respectively to the Eqns.~(2.1) and its combination with (2.2) in the above.
The condition (C) was first introduced, along with its first consequences, in \citeN[Section 4]{CrepeyElie16}:
it is explored systematically in this work. Whereas \citeN{CrepeySong15c} was in quest of generality and ``minimal conditions'', this paper aims at identifying a ``comfort zone'', provided by the condition (C) as the paper demonstrates, where all the standard apparatus of semimartingale calculus is equivalently available in both a larger stochastic basis $(\gg,\mathbb{Q})$ and a smaller one (making $\tau$ a stopping time)  $(\ff,\mathbb{P})$. But many problems reformulated under the reduced basis are simpler than their original formulation under the full stochastic basis, whence the benefit of the approach.

Specifically, assuming the condition (C) of an invariance time $\tau$ endowed with an intensity and a positive Az\'ema supermartingale as detailed in Section \ref{ss:C},
the present paper establishes a dictionary 
of transfer properties 
between the semimartingale calculi in the original and changed stochastic bases $(\gg,\mathbb{Q})$ and $(\ff,\mathbb{P})$: 
\begin{itemize}
\item Theorem \ref{t:exptransf} extends the classical reduced-form credit risk pricing formulas beyond the basic progressive enlargement of filtration setup where the Az\'ema supermartingale 
of $\tau$ has no martingale component;  
\item Theorem \ref{invarTH} establishes a bijection between the $(\gg,\mathbb{Q})$ (resp.~$(\gg,\mathbb{Q})$ continuous / $(\gg,\mathbb{Q})$ purely discontinuous) local martingales stopped before $\tau$ and the  $(\ff,\mathbb{P})$ (resp.~$(\ff,\mathbb{P})$ continuous / $(\ff,\mathbb{P})$ purely discontinuous)  local martingales;
\item Theorem \ref{itf-usualmartingale}
establishes the connection between 
stochastic integrals in the sense of local martingales in $(\gg,\mathbb{Q})$ and in $(\ff,\mathbb{P})$;
\item Theorem \ref{itf-randommeasure} establishes the connection between the $(\gg,\mathbb{Q})$ and $(\ff,\mathbb{P})$ random measures stochastic integrals;
\item Theorem \ref{mrtcoroll} establishes the correspondence between (weak or strong) $(\gg,\mathbb{Q})$  and  $(\ff,\mathbb{P})$ martingale representation properties;
\item Theorem \ref{tripletreduction} yields the relationship
between the $(\gg,\mathbb{Q})$ local characteristics of a $\gg$ semimartingale $X$ stopped before $\tau$
and the $(\ff,\mathbb{P})$ local characteristics of the $\ff$ semimartingale $X'$, called reduction of $X$, that coincides with $X$ before $\tau$; 
\item Theorems 
\ref{th:strm} and \ref{th:gene} state conditions under which Markov properties, transition semigroups and infinitesimal generators can be transferred between $(\gg,\mathbb{Q})$ and $(\ff,\mathbb{P})$; 
\item Theorems \ref{t:bsdeequiv} and \ref{c:eqp} show the equivalences, within various spaces of solutions, between a
nonstandard  $(\gg,\Q)$ backward SDE (BSDE)
stopped before $\tau$ and a reduced $(\ff,\P)$ BSDE with null terminal condition. 
\end{itemize} 
As illustrated in Section \ref{ss:ex},
the notion of invariance time is also related to various approaches that were introduced in the mathematical finance literature for coping with defaultable cash flows based on default intensities. These different approaches
could perhaps be related via the generalized Girsanov formulas of
\citeN{Kunita76} and \citeN{Yoeurp1985}.
We leave this for further research.

\appendix

	\section{Intensity Based Pricing Formulas, Survival Measure and Invariance Times}\label{ss:ex}

This section puts Theorem \ref{t:exptransf} (specifically, the formula \qr{t:etfp})  in perspective  with  \shortciteN[Proposition 1]{DuffieSchroderSkiadas96} 
and \citeN[Theorem 1]{CollinDufresneGoldsteinHugonnier2004}. 
This is done in the setup of the Markov model corresponding to the univariate case in Example \ref{e:dc}(ii), where the issues at stake can be understood based on Feynman-Kac representations. See \citeN{jeanblanc2020characteristics} for other renewed views on the seminal formulas of 
\shortciteN{DuffieSchroderSkiadas96}.

\subsection{The Univariate Dynamic Gaussian Copula Model}

We consider the following single-name version of the 
dynamic Gaussian copula model of
\citeN{CrepeySong16b}.
Let \beql{e:tau}\tau=\Psi\big(\int_0^{+\infty}\varsigma(s)dB_s
\big),\eeql 
where $\Psi$ is a continuously differentiable
increasing function from $\mathbb{R}$ to $(0,+\infty),$ 
$\varsigma$ is a Borel function on $\R_+$ such that 
$\int_0^{+\infty} \varsigma^2(u)du=1$, and $B$ is an $(\ff,\Q)$ Brownian motion, with $\ff$ taken as the \b{augmented natural} filtration of $B$.
The full
model filtration $\gg$ is given as the augmented filtration of the 
progressive enlargement of $\ff$ by
$\tau$. Hence, the random time $\tau$ is an $\cF_{\infty}$ measurable $\gg$ stopping time. 

By Theorem 2.2, Lemma 3.2, and Remark 4.1 in \citeN{CrepeySong16b}, the condition (C) holds in this setup,
for some probability measure  
\index{p@$\mathbb{P}$}$\mathbb{P}$ distinct from $\Q$ but equivalent to it on {${\mathfrak{F}}_{T}$}, on which $\mathbb{P}$ is uniquely determined through \qr{e:conda}. 
Let 
\beql{e:mk}
h_t=\ind_{\{t\geq \tau\}}\sp\index{m@$m$}m_t=\int_0^t\varsigma(s)dB_s\sp
\index{k@$k$}{k}_t= 
(h_t,\tau \wedge t )  \sp 
\mynu^2(t)= {\int_t^{+\infty} \varsigma ^2(s)ds},
\eeql 
and assume
$ \mynu$ positive for all $t$. \b{By application of results in \citeN{EKu}, one can show} that
the process $( {m} , {k} )$ is $(\gg,\Q)$ Markov. 
\b{\brem The reason why we introduce ($\tau \wedge t $) on top of the indicator $h_t$ in
${k}_t  $ is because of a dependence of the \b{post-$\tau$ behavior} of the model
on the value of $\tau$ itself. The state augmentation by ($\tau \wedge t$)
takes care of this path-dependence.~\finproof\erem}

By definition of $\mathtt{S}$ and  \qr{e:tau} of $\tau$, we have
\beql{e:imm}
\mathtt{S}_t=\Q(\tau
>t\,|\,\cF_t)=\Phi \big(\frac{\Psi^{-1}(t)-m_t}{\mynu (t)}\big)\sp t\in\R_+,\eeql
where $\Phi$ denotes the standard normal cdf.
The process on the right hand side of \qr{e:imm} has infinite variation.
This shows that
 the reference filtration
$\ff$ is not immersed into the full model filtration
$\gg$. This lack of immersion makes it more interesting from the point of view of the different approaches that we want to compare. 
This is our motivation for working in this particular model in this part.

Theorems 2.2 and 2.4 in
\citeN{CrepeySong16b}  
show the existence of
processes of the form
\beql{e:intensi}
\themu _t =\themu  (t, 
{m}_{t} ,k_t
) \mbox{ and }
\gamma _t = \gamma \, (t, 
{m}_{t} ,k_t
)  =\gamma _t \ind_{(0,\tau]}\sp t\in\R_+,\eeql
for 
continuous 
functions $\themu$ and $\gamma$ with linear growth in ${m}$, such that
\beql{e:wj}&\mbox{\index{w@$W$}$dW_t=dB_t- \themu_t dt$
is a $(\gg,\mathbb{Q})$ Brownian motion and}\\&
\mbox{the process $\gamma$
 is 
the $(\gg,\Q)$ intensity of $\tau$.}\eeql  

\noindent
The proof of the following result
is deferred to \sr{ss:theproof}.
\bp\label{p:codu} Let a process $m^{\star}$ satisfy
\beql{e:mst}dm^{\star}_t=\varsigma (t) \big(dW^{\star}_t + {\beta} (t,m^{\star}_t, (0,t)  ) dt\big)\sp 0\le t\le T,\eeql
starting from $m^{\star}_0=0$, 
for some Brownian motion $W^{\star}$ with respect to some stochastic basis $(\gg^{\star},\Q^{\star}).$
Denoting the $\Q^{\star}$ expectation by ${\E}^{\star}$, we have, for any %
bounded Borel function $G(t,m),$  
\beql{e:vDuff}
&\E\big[\ind_{\{\tau\leq T\}} G(\tau,m_\tau )  \big]={\E}^{\star} \Big[\int_0^T e^{-\int_0^t  \gamma(s,m^{\star}_s, (0,s))  ds  }\gamma(t,m^{\star}_t, (0,t))  G(t,m^{\star}_t) dt  \Big].
\eeql
\ep

\subsection{Discussion}
From \qr{e:mk} and \qr{e:intensi}-\qr{e:wj}, it holds that
\beql{e:m}dm_t=\varsigma (t) \big(dW_t + {\beta} (t,m_t, k_t  )dt\big)\sp t\in\R_+,\eeql
which, for $t\ge\tau$ (so that $k_t=(1,\tau)$), diverges from the specification \qr{e:mst}. Hence the quadruplet
$(m,W,\gg,\Q)$ is not an eligible choice for $(m^\star,W^\star,\gg^\star,\Q^\star)$ in Proposition \ref{p:codu}. As a consequence,
we expect a contrario from Proposition \ref{p:codu} that
\beql{e:dufftocorr}\E\big[\ind_{\{\tau<T\}} G(\tau,m_\tau )  \big]\neq {\E} \Big[\int_0^{T} e^{-\int_0^t  \gamma(s,m_s,(0,s))  ds  }\gamma(t,m_t,(0,t))  G(t,m_t) dt  \Big]  \eeql
(except in special cases, including obviously $G=0$).
In fact, let $$
V_t= {\E}  \Big[\int_t^{T} e^{-\int_t^s  \gamma(u,m_u, (0,u))  du  }\gamma(s,m_s, (0,s))  G(s,m_s) ds\,\Big|\,\cG_t \Big]\sp t\in\R_+,$$ 
so that $V_0$ is equal to the right hand side in \qr{e:dufftocorr}.
By an application of \shortciteN[Proposition 1]{DuffieSchroderSkiadas96} with $X=r=0$
and $h_{\cdot}=\gamma(\cdot,m_{\cdot}, (0,\cdot))$ on $[0,T]$ there (noting that any process coinciding with the $(\gg,\Q)$ intensity of $\tau$ before $\tau$ can be used as a process $h$ in their setup), we have
\beql{e:df}\E\big[\ind_{\{\tau<T\}} G(\tau,m_\tau )  \big]=V_0-{\E} (V_{\tau}  -V_{\tau-}) .\eeql
In a basic immersed setup, $\E(V_{\tau}  -V_{\tau-})$ vanishes and equality actually holds in \qr{e:dufftocorr} for any $G$: see 
the comments before Section 3 in \shortciteN{DuffieSchroderSkiadas96},
page 1379 in \shortciteN{CollinDufresneGoldsteinHugonnier2004},
or following (3.22), (H.3) and Proposition 6.1 in \citeN{BieleckiRutkowski00}).
But beyond this immersion case, ${\E} (V_{\tau}  -V_{\tau-})$
is typically nonnull and intractable, whence inequality in \eqref{e:dufftocorr} (except in special cases such as $G=0$).

Instead, an eligible choice in Proposition \ref{p:codu}
consists in using $m^\sta=m$ ,
$W^\sta=W$ as per \eqref{e:wj}, $\mathbb{Q}^\sta=$ the so called  survival measure $\mathbb{S}$
 with  $(\gg,\mathbb{Q})$ density process 
$$e^{\int_0^{\cdot\wedge T} \gamma(u,m_u, (0,u))  du }\indi{\tau>\cdot\wedge T}$$ (assuming $e^{\int_0^{\tau\wedge T} \gamma(u,m_u, (0,u))  du }
$ integrable under $\Q$), and
$\gg^\sta=$ the 
 $\mathbb{S}$ augmentation $\bar{\gg}$ of $\gg$, obtained  
by adding \b{to each $\cG_t$} all the $\mathbb{S}$ null sets $A\in\cA$ such that $A\subseteq  \{\tau\leq T\}$.
Indeed, as noted in \shortciteN[Lemma 1(i)]{CollinDufresneGoldsteinHugonnier2004},
$\ind_{[\tau,+\infty)}=0$ holds $\mathbb{S}$ almost surely on $[0,T]$, hence
\eqref{e:mst}
holds  in this setup, while \shortciteN[Lemma 1(ii)]{CollinDufresneGoldsteinHugonnier2004} shows that 
$W^\sta=W$ is a $(\gg^{\star},\Q^{\star})=(\bar{\gg},\mathbb{S})$ Brownian motion.
The corresponding specification of the formula \eqref{e:vDuff} corresponds to \shortciteN[Theorem 1]{CollinDufresneGoldsteinHugonnier2004}.
This ``survival measure'' idea and terminology were first introduced in 
\citeN{Schoenbucher99,Schoenbucher04}.
One can thus fix the discrepancy in \qr{e:dufftocorr} (in a progressive enlargement of filtration setup without immersion) by singularly changing the probability measure $\mathbb{Q}$ to $\mathbb{Q}^\sta=\mathbb{S}$, while sticking to the original model filtration $\gg$ (or, more precisely, resorting to its
 $\mathbb{S}$ augmentation $\bar{\gg}$).

Another eligible choice for $(m^\star,W^\star,\gg^\star,\Q^\star)$ in Proposition \ref{p:codu}
is
$$m^\sta=m\sp  \index{w@$W$} d{W}^\sta_t=dB_t- \themu(t,m_t,(0,t)) dt\sp \gg^\sta=\ff\sp\mathbb{Q}^\sta=\P.$$ 
Indeed, as it follows from
Lemma 3.5 and Section 4.4 in
\citeN{CrepeySong16b},
this process ${W}^\sta$
 is an $( \ff, \P)$ Brownian motion.
 Hence, just like $(m,
B-  \themu(\cdot,m_\cdot,k_\cdot)\centerdot\boldsymbol\lambda,
\bar{\gg},\mathbb{S})$  in the previous specification,
 $(m,B- \themu(\cdot,m_\cdot,(0,\cdot) )
\centerdot\boldsymbol\lambda , {\ff} ,\P)$ satisfies all the conditions required on $(m^\sta,W^\sta, {\gg}^\sta,\Q^\sta)$ in Proposition \ref{p:codu}.
The formula \qr{e:vDuff} corresponding to this second legit specification for  $(m^\sta,W^\sta, {\gg}^\sta,\Q^\sta)$
is none other than our 
formula \qr{t:etfp} (for $t=0$). 

The approach $(m^\sta,W^\sta, {\gg}^\sta,\Q^\sta)=(m,B- \themu(\cdot,m_\cdot,(0,\cdot) )
\centerdot\boldsymbol\lambda , {\ff} ,\P)$ of this paper thus fixes the discrepancy in \qr{e:dufftocorr} (in a non-immersed setup)  by reducing the filtration from $\gg$ to a smaller $\ff,$ while changing the probability measure ``as little as possible'', i.e.~equivalently on $\cF_T$.
\shortciteN{CollinDufresneGoldsteinHugonnier2004}'s  choice $(m^\sta,W^\sta, {\gg}^\sta,\Q^\sta)=(m,
B-  \themu(\cdot,m_\cdot,k_\cdot)\centerdot\boldsymbol\lambda,
\bar{\gg},\mathbb{S})$ does the opposite, touching the filtration as little as possible but changing the measure singularly
(in a basic
immersive setup, an invariance time approach would not change $\Q$ at all, whereas \shortciteN{CollinDufresneGoldsteinHugonnier2004}'s measure change would still be singular).
But \shortciteN{CollinDufresneGoldsteinHugonnier2004}
only provide
a transfer of conditional expectation formulas.
The present paper demonstrates how the invariance times approach, instead,
results in a transfer of semimartingale calculus as a whole.
One concrete motivation 
for this work
is the solution of BSDEs stopped before their terminal time.
As Section \ref{s:bsde} illustrates, 
in order to deal with these,
conditional expectation formulas are not enough:
the entire semimartingale calculus of this paper is required.

\subsection{Proof of Proposition \ref{p:codu}\label{ss:theproof}}
By the $(\gg,\Q)$ Markov property of the process $( {m} , {k} )$,
\b{noting that $\tau$ is the hitting time of 1 by the $h$ component of the process $k$},
we have
\beql{e:G}
&\E\big[\ind_{\{\tau\le T\}} G(\tau,m_\tau ) \,|\,\cG_t\big]=\E\big[\ind_{\{\tau<T\}} G(\tau,m_\tau ) \,|\,(m_t ,k_t ) \big]=\\&\qqq v(t,m_t,k_t)=u_{h_t}(t,m_t)\sp t \in [0,\tau \wedge T],\eeql
for suitable Borel bounded functions $v(t,m,k)$ and $u_{h}(t,m)=v(t,m,k=(h,t))$.
As a $(\gg,\Q)$ martingale, the process $u_{h_t}(t,m_t)$, $t \in [0, \tau \wedge T]$, has a vanishing $(\gg,\Q)$ drift.
Hence, by an application of the It\^o formula to this process, using \qr{e:wj}, the pair function $u=(u_0(t,m),u_1(t,m))$
formally solves $u_1  =G, u_0(T,\cdot )= 0$ and 
\beql{etheVI0remmarksbar} 
&\partial_t u_0(t,m ) + \varsigma (t)   {\beta} (t,m , (0,t))
\partial_{m  } u_0(t,m )  +\frac{\varsigma (t) ^2}{2} 
\partial^2_{m^2} u_0(t,m )\\&\qqq  + {\gamma}(t,m,0 )\big[G(t,m ) -u_0 (t,m )\big]
 = 0 \sp t<T,  m  \in \R.
\eeql
 \brem At least, the above holds
assuming $u$ regular enough for applicability of the It\^o formula.
Given the discussion nature of this appendix, we content ourselves with the above formal argument, 
without introducing weak solutions of \eqref{etheVI0remmarksbar}.
\erem
Putting together \qr{e:G} and the Feynman-Kac representation of the solution $u_0$ of \qr{etheVI0remmarksbar} at the origin yields
\bel
&\E\big[\ind_{\{\tau\le T\}} G(\tau,m_\tau )  \big]=
u_0(0,0)\\&\qqq
={\E}^{\star} \Big[\int_0^T e^{-\int_0^t  \gamma(s,m^{\star}_s, (0,s))  ds  }\gamma(t,m^{\star}_t, (0,t))  G(t,m^{\star}_t) dt  \Big], 
\eel
for any process $m^{\star}$ as stated in the proposition, which is therefore proven.

\end{document}